\documentclass[journal]{IEEEtran}
\usepackage{indentfirst}
\usepackage{mathrsfs}
\usepackage{graphicx}
\usepackage{epstopdf}
\usepackage{amsmath, amssymb, amsthm, bm}
\usepackage{leftidx}
\usepackage{extarrows}
\usepackage{stfloats}
\usepackage{picinpar}
\usepackage{enumerate}
\usepackage{epsfig}
\usepackage{latexsym}
\usepackage{amsfonts}
\usepackage{enumerate}
\usepackage{graphics}
\usepackage{float}
\usepackage{pict2e}
\usepackage{tikz}
\usepackage{lineno}
\usepackage[pdftex,colorlinks,linkcolor=blue,anchorcolor=blue,citecolor=blue]{hyperref}
\usepackage{moreverb}
\usepackage{color}
\usepackage[nosort]{cite}
\usepackage{subfigure}
\usepackage{multirow}
\usetikzlibrary{arrows,automata}
\usepackage{colortbl}
\usepackage{wrapfig}
\usepackage{threeparttable}
\usepackage{cuted}
\usepackage{units}
\usepackage{bbm}
\usepackage{bbding}
\usetikzlibrary{arrows,automata}

\usepackage[linesnumbered,ruled,vlined]{algorithm2e}

\definecolor{myGray}{rgb}{.95,.95,.95}

\newfont{\bb}{msbm10}

\newcommand{\tr}{^{\sf T}}
\newcommand{\M}[1]{{\bf{#1}}}

\newcommand{\m}[1]{{\mathrm{#1}}}

\newtheorem{thm}{Theorem}
\newtheorem{ass}{Assumption}

\newtheorem{lem}{Lemma}

\newtheorem{prop}{Proposition}

\ifCLASSINFOpdf

\else

\fi

\newcommand{\Exp}[1]{\mathbb{E}\!\left[ #1 \right]}

\begin{document}
\title{A Proximal Gradient Method With Probabilistic Multi-Gossip Communications for Decentralized Composite Optimization}

\author{Luyao Guo,
Luqing Wang,
Xinli Shi, \IEEEmembership{Senior Member, IEEE},
Jinde Cao, \IEEEmembership{Fellow, IEEE}
\thanks{This work was supported in part by the National Key R\&D Project of China  under Grant Nos. 2020YFA0714300, and the National Natural Science Foundation of China under Grant Nos. 62576098 and 62473098. (Corresponding author: Jinde Cao.)}

\thanks{Luyao Guo is with the School of Computer Science and Engineering, Suzhou University of Technology, Suzhou 215500, China, and also with the School of Mathematics, Southeast University, Nanjing 210096, China
(e-mail: \href{mailto:ly_guo@seu.edu.cn}{{ly\_guo@seu.edu.cn}}).}
\thanks{Luqing Wang is with the School of Computer Science \& Engineering, Southeast University, Nanjing 210096, China (e-mail: \href{mailto:luqing_wang@seu.edu.cn}{{luqing\_wang@seu.edu.cn}}).}
\thanks{Xinli Shi is with the School of Cyber Science and Engineering, Southeast University, Nanjing 210096, China
(e-mail: \href{mailto:xinli_shi@seu.edu.cn}{{xinli\_shi@seu.edu.cn}}).}
\thanks{Jinde Cao is with the School of Mathematics, Southeast University, Nanjing 210096, China, and also with the Purple Mountain Laboratories, Nanjing 211111, China. (e-mail: \href{mailto:jdcao@seu.edu.cn}{{jdcao@seu.edu.cn}}).}
}

\date{}

\maketitle

\begin{abstract}
Decentralized optimization methods with local updates have recently gained attention for their provable ability to communication acceleration. In these methods, nodes perform several iterations of local computations between the communication rounds. Nevertheless, this capability is effective only when the network is sufficiently well-connected and the loss function is smooth. In this paper, we propose a communication-efficient method $\textsc{MG-Skip}$ with probabilistic local updates and multi-gossip communications for decentralized composite (smooth + nonsmooth) optimization, whose stepsize is independent of the number of local updates and the network topology. For any undirected and connected networks, $\textsc{MG-Skip}$ allows for the multi-gossip communications to be skipped in most iterations in the strongly convex setting, while its computation complexity is $\mathcal{O}(\kappa \log \nicefrac{1}{\epsilon})$ and communication complexity is only $\mathcal{O}(\sqrt{\nicefrac{\kappa}{(1-\rho)}} \log \nicefrac{1}{\epsilon})$, where $\kappa$ is the condition number of the loss function, $\rho$ reflects the connectivity of the network topology, and $\epsilon$ is the target accuracy. The theoretical results indicate that $\textsc{MG-Skip}$ achieves provable communication acceleration, thereby validating the advantages of local updates in the nonsmooth setting.
\end{abstract}


\section{Introduction}
Consider the composite optimization problem over undirected and connected networks with $n$ agents
\begin{equation}\label{EQ:Problem1}
\begin{aligned}
\m{x}^\star   &= \arg\min_{\m{x}\in \mathbb{R}^d}\left\{h(\m{x}):=\frac{1}{n}\sum_{i=1}^{n}f_i(\m{x})+r(\m{x})\right\},
\end{aligned}
\end{equation}
where $f_i:\mathbb{R}^d\rightarrow \mathbb{R}$ is the local loss function accessed by the node $i$, and $r:\mathbb{R}^d\rightarrow \mathbb{R}\cup\{\infty\}$ is a global loss function. In this setup, a network of nodes (also referred to as agents, workers, or clients) collaboratively seeks to minimize the average of the nodes' objectives. We give the following assumption throughout this paper.
\begin{ass}\label{ASS1}
We assume that
\begin{enumerate}
  \item $f_i$ is $\mu_i$-strongly convex and $L_i$-smooth (with $\mu_i>0$ and $L_i\geq0$). Define $L:=\max_i\{L_i\}$, $\mu:=\min_i\{\mu_i\}$, and the condition value of $f_i$ as $\kappa:= \nicefrac{L}{\mu}$.
  \item $r$ is convex closed and proper. Moreover, $r$ is proximable, i.e., with proper $\alpha>0$, its proximal mapping
  $$
  \mathrm{prox}_{\alpha r}(\m{y})=\mathrm{argmin}_{\m{x}\in\mathbb{R}^d}\left\{\alpha r(\m{x})+\nicefrac{1}{2}\|\m{x}-\m{y}\|^2\right\},
  $$
has an analytical solution or can be computed efficiently.
\end{enumerate}
\end{ass}

\begin{figure}[!t]
  \centering
  \setlength{\abovecaptionskip}{-2pt}
  \includegraphics[width=1\linewidth]{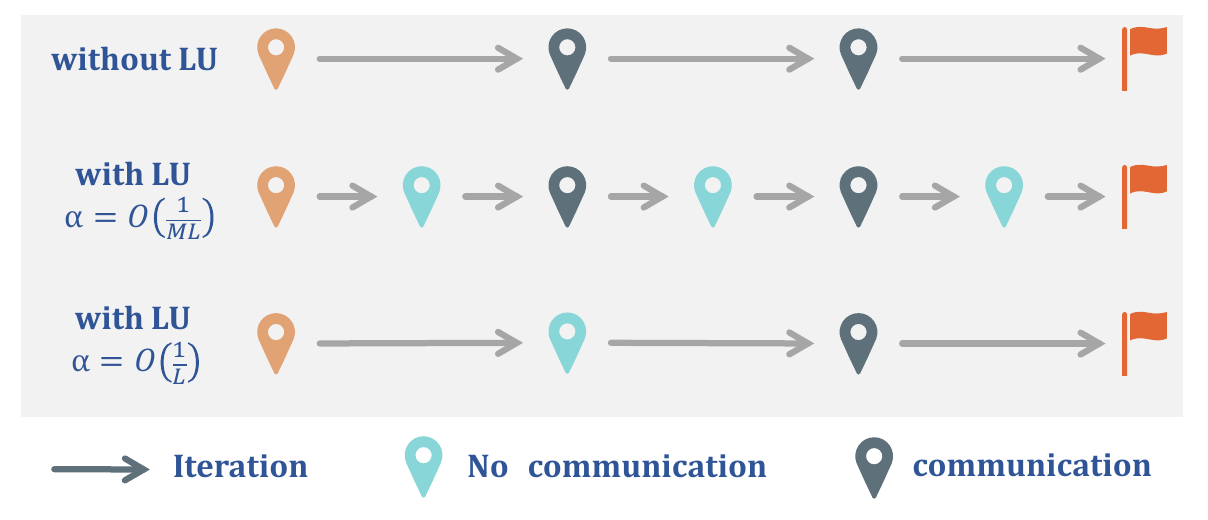}
  \caption{Sketch of the execution of decentralized optimization algorithms.}
  \label{Sketch-step}
\end{figure}

This setting is general and is encountered in various application fields, including multi-agent control, wireless communication, and distributed machine learning \cite{Nedic20151,Nedic2018,Jakovetic2019,DADMM}. Various algorithms have been proposed to solve problem \eqref{EQ:Problem1} in a decentralized manner, such as EXTRA \cite{EXTRA,PGEXTRA}, DIGing \cite{DIGing}, GT \cite{Harnessing}, NEXT \cite{NEXT}, AccNEXT \cite{AccNEXT}, NIDS \cite{NIDS} (or ED \cite{ExactDiffusion}), SONATA \cite{SONATA}, PMGT-VR \cite{PMGT-VR}, D-iPGM \cite{Guo2022}, DISA \cite{Guo2023}, and VR-DPPD \cite{LU2022}. To retain the linear convergence rate for these algorithms to problem \eqref{EQ:Problem1}, \cite{Sulaiman2021,Xu2021} propose their proximal gradient extension. These methods do not rely on a central coordinator and that communicate only with neighbors in an arbitrary communication topology. Nevertheless, decentralized optimization algorithms may still face challenges arising from communication bottlenecks.

\begin{table*}[!t]
\renewcommand\arraystretch{2}
\begin{center}
\caption{A comparison with existing decentralized optimization methods. $T$ and $C$ denote computation and communication complexity, respectively. $(1-\rho)$-InS denotes $(1-\rho)$-independent stepsize. $M$-InS denotes $M$-independent stepsize. LU = local updates.}
\begin{threeparttable}
\scalebox{0.92}{
\begin{tabular}{ccccccccc}
\hline
Algorithm &Stepsize & $T$ & $C$ &$r\neq0$&$(1-\rho)$-InS&$M$-InS&LU&opt.? \tnote{a}  \\
\hline
Prox-EXTRA \cite{EXTRA,Sulaiman2021,Xu2021} &$\mathcal{O}\left(\frac{1-\rho}{ L}\right)$&$\mathcal{O}\left(\frac{\kappa}{1-\rho}\log \nicefrac{1}{\epsilon}\right)$&$\mathcal{O}\left({\frac{\kappa}{1-\rho}}\log \nicefrac{1}{\epsilon}\right)$&\textcolor[rgb]{0,0.6,0}{\Checkmark}&\textcolor[rgb]{0.7,0,0}{\XSolidBrush}&\textcolor[rgb]{0.7,0,0}{\XSolidBrush}&-&\textcolor[rgb]{0.7,0,0}{\XSolidBrush}\\
Prox-GT \cite{Harnessing,Sulaiman2021,Xu2021} &$\mathcal{O}\left(\frac{(1-\rho)^2}{ L}\right)$&$\mathcal{O}\left(\frac{\kappa}{(1-\rho)^2}\log \nicefrac{1}{\epsilon}\right)$&$\mathcal{O}\left(\frac{\kappa}{(1-\rho)^2}\log \nicefrac{1}{\epsilon}\right)$&\textcolor[rgb]{0,0.6,0}{\Checkmark}&\textcolor[rgb]{0.7,0,0}{\XSolidBrush}&\textcolor[rgb]{0.7,0,0}{\XSolidBrush}&-&\textcolor[rgb]{0.7,0,0}{\XSolidBrush}\\
Prox-DIGing \cite{DIGing,Sulaiman2021,Xu2021} &$\mathcal{O}\left(\frac{(1-\rho)^2}{ L}\right)$&$\mathcal{O}\left(\frac{\kappa}{(1-\rho)^2}\log \nicefrac{1}{\epsilon}\right)$&$\mathcal{O}\left(\frac{\kappa}{(1-\rho)^2}\log \nicefrac{1}{\epsilon}\right)$&\textcolor[rgb]{0,0.6,0}{\Checkmark}&\textcolor[rgb]{0.7,0,0}{\XSolidBrush}&\textcolor[rgb]{0.7,0,0}{\XSolidBrush}&-&\textcolor[rgb]{0.7,0,0}{\XSolidBrush}\\
Prox-NIDS \cite{NIDS,Sulaiman2021,Xu2021} & $\mathcal{O}\left(\frac{1}{L}\right)$&$\mathcal{O}\left(\frac{\kappa}{1-\rho}\log \nicefrac{1}{\epsilon}\right)$&$\mathcal{O}\left(\frac{\kappa}{1-\rho}\log \nicefrac{1}{\epsilon}\right)$&\textcolor[rgb]{0,0.6,0}{\Checkmark}&\textcolor[rgb]{0,0.6,0}{\Checkmark}&\textcolor[rgb]{0.7,0,0}{\XSolidBrush}&-&\textcolor[rgb]{0.7,0,0}{\XSolidBrush}\\
MG-SONATA \cite{SONATA} / PMGT-VR \cite{PMGT-VR}  &$\mathcal{O}\left(\frac{1}{L}\right)$&$\mathcal{O}(\kappa\log \nicefrac{1}{\epsilon})$&$\mathcal{O}\left(\frac{\kappa}{\sqrt{1-\rho}}\log \nicefrac{1}{\epsilon}\right)$&\textcolor[rgb]{0,0.6,0}{\Checkmark}&\textcolor[rgb]{0,0.6,0}{\Checkmark}&\textcolor[rgb]{0.7,0,0}{\XSolidBrush}&-&\textcolor[rgb]{0.7,0,0}{\XSolidBrush}\\
\hline
local-GT \cite{Liu2023}&$\mathcal{O}\left(\frac{(1-\rho)^2}{ML}\right)$&$\mathcal{O}\left(\frac{M\kappa}{(1-\rho)^2}\log \nicefrac{1}{\epsilon}\right)$&$O\left(\frac{\kappa}{(1-\rho)^2\log \nicefrac{1}{\epsilon}}\right)$&\textcolor[rgb]{0.7,0,0}{\XSolidBrush}&\textcolor[rgb]{0.7,0,0}{\XSolidBrush}&\textcolor[rgb]{0.7,0,0}{\XSolidBrush}&\textcolor[rgb]{0,0.6,0}{\Checkmark}&\textcolor[rgb]{0.7,0,0}{\XSolidBrush}\\
LED \cite{Alghunaim2023}&$\mathcal{O}\left(\frac{(1-\rho)\mu}{ML^2}\right)$&$\mathcal{O}\left(\frac{M\kappa^2}{1-\rho}\log \nicefrac{1}{\epsilon}\right)$&$\mathcal{O}\left(\frac{\kappa^2}{1-\rho}\log \nicefrac{1}{\epsilon}\right)$&\textcolor[rgb]{0.7,0,0}{\XSolidBrush}&\textcolor[rgb]{0.7,0,0}{\XSolidBrush}&\textcolor[rgb]{0.7,0,0}{\XSolidBrush}&\textcolor[rgb]{0,0.6,0}{\Checkmark}&\textcolor[rgb]{0.7,0,0}{\XSolidBrush}\\
ProxSkip \cite{ProxSkip}&${\mathcal{O}}\left(\frac{1}{L}\right)$&$\mathcal{O}\left(\frac{\kappa}{1-\rho}\log \nicefrac{1}{\epsilon}\right)$&$\mathcal{O}\left(\sqrt{\frac{\kappa}{1-\rho}}\log \nicefrac{1}{\epsilon}\right)$ \tnote{b}&\textcolor[rgb]{0.7,0,0}{\XSolidBrush}&\textcolor[rgb]{0,0.6,0}{\Checkmark}&\textcolor[rgb]{0,0.6,0}{\Checkmark}&\textcolor[rgb]{0,0.6,0}{\Checkmark}&$\rho\leq 1-\nicefrac{1}{\kappa}$\\
\hline
\rowcolor[gray]{0.9} $\textsc{MG-Skip}$&${\mathcal{O}}\left(\frac{1}{L}\right)$&$\mathcal{O}\left(\kappa\log \nicefrac{1}{\epsilon}\right)$&$\mathcal{O}\left(\sqrt{\frac{\kappa}{1-\rho}}\log \nicefrac{1}{\epsilon}\right)$&\textcolor[rgb]{0,0.6,0}{\Checkmark}&\textcolor[rgb]{0,0.6,0}{\Checkmark}&\textcolor[rgb]{0,0.6,0}{\Checkmark}&\textcolor[rgb]{0,0.6,0}{\Checkmark}&\textcolor[rgb]{0,0.6,0}{\Checkmark}\\
\hline
\end{tabular}}
\begin{tablenotes}
  \footnotesize
   \item[a] Is the optimal communication complexity established in \cite{Scaman2017} achieved?
   \item[b] The communication complexity holds only for $\rho\leq 1-\nicefrac{1}{\kappa}$.
\end{tablenotes}
 \end{threeparttable}
 \label{TableCOMPA}
\end{center}
\end{table*}

To reduce communication costs, the recent surge in local updates methods has ignited considerable interest \cite{Liu2023,Alghunaim2023,ProxSkip}. In these algorithms, a collection of nodes perform several iterations of local computations between the communication rounds. In \cite{Liu2023} and \cite{Alghunaim2023}, the local-GT and LED algorithms have been respectively introduced, both of which integrate deterministic periodic local updates. Although these algorithms demonstrate the advantages of employing local updates in stochastic settings, they do not offer a provable benefit in deterministic settings. The main reason for this phenomenon is that the allowable stepsizes for both local-GT \cite{Liu2023} and LED \cite{Alghunaim2023} are of the order $\mathcal{O}\left(\frac{1}{M L}\right)$, where $M$ represents the number of local updates, indicating that more local updates lead to smaller stepsizes (see Fig \ref{Sketch-step}). In contrast, ProxSkip \cite{ProxSkip}, which incorporates probabilistic local updates, can accelerate communication and achieve the optimal communication complexity of $\mathcal{O}(\sqrt{\nicefrac{\kappa}{(1-\rho)}}\mathrm{log}\nicefrac{1}{\epsilon})$ \cite{Scaman2017} in deterministic settings, where $\rho$ reflects the connectivity of the network topology; specially, if $\rho=0$ the network is fully connected. This achievement is achieved due to the allowable stepsize of ProxSkip is of the order $\mathcal{O}\left(\frac{1}{ L}\right)$, which is independent of the number of local updates, and it is linear convergent.

Nonetheless, it has been demonstrated both theoretically, as indicated in the work by \cite{ProxSkip}, and experimentally, as evidenced by \cite{Alghunaim2023}, that the existing decentralized methods with local updates \cite{Liu2023,Alghunaim2023,ProxSkip} can exclusively achieve communication acceleration when the network exhibits sufficient connectivity, specifically when $\rho\leq 1-\nicefrac{1}{\kappa}$. Furthermore, it is worth noting that these methods \cite{Liu2023,Alghunaim2023,ProxSkip} are only applicable under the smooth case where $r(\m{x})\equiv 0$. Thus, it naturally raises a pivotal inquiry: for decentralized composite optimization \eqref{EQ:Problem1} on general undirected and connected networks, is it possible to achieve provable communication acceleration using probabilistic local updates?

In this paper, we answer this question in the affirmative and propose a novel proximal gradient decentralized method with random \underline{M}ulti-\underline{G}ossip communication \underline{Skip}pings ($\textsc{MG-Skip}$) for the nonsmooth decentralized optimization problem \eqref{EQ:Problem1}.
We conduct a comparison of $\textsc{MG-Skip}$ with existing methods in Table \ref{TableCOMPA}, where MG-SONATA \cite{SONATA} denotes the SONATA algorithm with extra averaging steps by Chebyshev acceleration. The main contributions of this paper are outlined below.
\begin{itemize}
  \item We introduce a novel algorithm named by $\textsc{MG-Skip}$ for problem \eqref{EQ:Problem1}, where communication occurs with a probability of $p\in(0,1]$. With network-independent and local updates number-independent stepsizes, we prove the linear convergence of $\textsc{MG-Skip}$. To the best of our knowledge, $\textsc{MG-Skip}$ is the first linear convergent decentralized method with probabilistic local updates for nonsmooth settings.
  \item We prove that $\textsc{MG-Skip}$ not only inherits the advantages of ProxSkip \cite{ProxSkip}---provable communication-acceleration and often skip communication for free,
  but also overcomes the dependence on the network---no requirements for well-connected network.
  \item When $p\in[\nicefrac{1}{\sqrt{\kappa}},1]$, the computation complexity $T$ and communication complexity $C$ of $\textsc{MG-Skip}$ are
    \begin{align*}
    &T=\mathcal{O}\left(\kappa\log \nicefrac{1}{\epsilon}\right),\quad C=\mathcal{O}\left(\frac{p\kappa}{\sqrt{1-\rho}}\log \nicefrac{1}{\epsilon}\right).
    \end{align*}
  It implies that $\textsc{MG-Skip}$ achieves provable communication acceleration. When setting $p=\nicefrac{1}{\sqrt{\kappa}}$, it achieves the optimal communication complexity \cite{Scaman2017}, i.e.,
  $$C=\mathcal{O}\left(\sqrt{\frac{\kappa}{1-\rho}}\log \nicefrac{1}{\epsilon}\right).$$
\end{itemize}

This paper is organized as follows. In Section \ref{SEC3}, we introduce $\textsc{MG-Skip}$ and discuss the relation to existing methods. In Section \ref{SEC4}, we give the main theoretical results. Finally, several numerical simulations are implemented in Section \ref{SEC5} and conclusions are given in Section \ref{SEC6}.

\section{The Proposed Algorithm}\label{SEC3}
For any $n\times m$ matrices $\mathbf{A}$ and $\mathbf{C}$, their inner product is denoted as $\langle\mathbf{A},\mathbf{C}\rangle=\mathrm{Trace}(\mathbf{A}\tr\mathbf{C})$, the Frobenius norm is given by $\|\mathbf{A}\|_\mathrm{F}$, while the spectral norm is given by $\|\mathbf{A}\|$. $\sigma_m(\mathbf{A})$ denotes the smallest none zero eigenvalue of $\mathbf{A}$ and $\varpi(\mathbf{A})$ denotes the spectral radius of matrix $\mathbf{A}$. Let $\mathrm{null}(\cdot)$ and $\mathrm{span}(\cdot)$ be the range space and null space of the argument vector/matrix, respectively. The $1$-norm and Euclidean norm of a vector are denoted as $\|\cdot\|_1$ and $\|\cdot\|$. $0$ and $\M{I}$ denote the null matrix and the identity matrix of appropriate dimensions, respectively. Denote $\M{1}$ as a vector with each component being one. $\lfloor\cdot\rfloor$ denotes round toward $-\infty$.
\subsection{Network Graph}
Let $\M{x}_i^t\in\mathbb{R}^d$ be the local state of node $i$ at the $t$-th iteration. Moreover, we define the following notations
\begin{align*}
\M{x}:=&\ \left[\m{x}_{1},\m{x}_{2},\cdots,\m{x}_{n}\right]\tr\in\mathbb{R}^{n\times d},\\
\M{x}^t:=&\ \left[\m{x}_{1}^t,\m{x}_{2}^t,\cdots,\m{x}_{n}^t\right]\tr\in\mathbb{R}^{n\times d},\\
F(\M{x}):=&\sum_{i=1}^n  f_i(\m{x}_i),\ R(\M{x}):=\sum_{i=1}^n r(\m{x}_i),\\
\nabla F(\M{x}^t):=&\ \left[\nabla f_1(\m{x}^t_1),\nabla f_2(\m{x}^t_2),\cdots,\nabla f_n(\m{x}^t_n)\right]\tr\in \mathbb{R}^{n\times d},
\end{align*}
which collect all local variables and gradients across the network. In this work, we focus on decentralized scenarios, where a group of $n$ nodes is interconnected by an undirected graph with a set of edges $\mathcal{E}$, where node $i$ is connected to node $j$ if $(i,j) \in \mathcal{E}$. Let $\mathcal{N}_i$ denote the set of neighbours of agent $i$ including itself. We introduce the global mixing matrix $\mathbf{W}=[W_{ij}]_{i,j=1}^n\in\mathbb{R}^{n\times n}$, where $W_{ij}=W_{ji}=0$ if $(i,j) \notin \mathcal{E}$, and $W_{ij}>0$ otherwise. We impose the following standard assumption on the weight matrix.
\begin{ass}\label{MixingMatrix}
The network is undirected and connected, and the weight matrix $\M{W}$ is symmetric and doubly stochastic, i.e., $\M{W}=\M{W}\tr$ and $\M{W1}=\M{1}$.
\end{ass}
Let $\lambda_1=1$ denote the largest eigenvalue of the mixing matrix $\mathbf{W}$, and the remaining eigenvalues are denoted as
$$
1>\lambda_2\geq\lambda_3\geq\cdots\geq\lambda_n>-1.
$$
If we let
$\rho=\varpi(\M{W}-\nicefrac{1}{n}\M{11}\tr)$,
where $\varpi(\M{W}-\nicefrac{1}{n}\M{11}\tr)$ denotes the spectral radius of matrix $\M{W}-\nicefrac{1}{n}\M{11}\tr$,  then from \cite[Remark 1]{Yuan2021}, we have
$$
\rho=\max\{|\lambda_2|,|\lambda_n|\}\in(0,1).
$$
The network quantity $\rho$ is called the spectral gap of $\mathbf{W}$, which reflects the connectivity of the network topology. Based on this weight matrix $\M{W}$, we design the following multi-round weight matrix $\bar{\M{M}}:=\M{M}_K$, where $K$ denotes the required consecutive rounds of communications. Initial
$\M{M}_{-1}=\M{M}_0=\M{I}$, generate $\M{M}_K\in\mathbb{R}^{n\times n}$ by:
\begin{align*}
\text{\textbf{ for }} k=0,1,\ldots,K-1,\ \M{M}_{k+1}=(1+\eta)\M{W}\M{M}_k-\eta\M{M}_{k-1},
\end{align*}
where $\eta=\nicefrac{\big(1-\sqrt{1-\rho^2}\big)}{\big(1+\sqrt{1+\rho^2}\big)}$. Then, we have
\begin{prop}\label{Pro1}
With Assumption \ref{MixingMatrix}, it holds that
\begin{itemize}
  \item $\bar{\M{M}}$ is symmetric and doubly stochastic, i.e., $\bar{\M{M}}=\bar{\M{M}}\tr$ and $\bar{\M{M}}\M{1}=\M{1}$, and $\varpi(\bar{\M{M}}-\nicefrac{1}{n}\M{11}\tr)\leq\sqrt{2}(1-\sqrt{1-\rho})^K$.
  \item Setting $K=\lfloor\nicefrac{1}{\sqrt{1-\rho}}\rfloor$, we have $\sigma_m(\M{I}-\bar{\M{M}})\geq\nicefrac{2}{5}$,
where $\sigma_m(\M{I}-\bar{\M{M}})$ is the smallest none zero eigenvalue of $\M{I}-\bar{\M{M}}$.
\end{itemize}
\end{prop}
\begin{IEEEproof}
See Appendix \ref{AP:Pro1}.
\end{IEEEproof}

\subsection{Algorithm Development and Intuition}
By Proposition \ref{Pro1}, it holds that
\begin{align*}
\nicefrac{1}{2}(\M{I}-\bar{\M{M}})\M{x}=0&\Longleftrightarrow \sqrt{\nicefrac{1}{2}(\M{I}-\bar{\M{M}})}\M{x}=0\\
&\Longleftrightarrow  \m{x}_1=\cdots=\m{x}_n.
\end{align*}
The problem \eqref{EQ:Problem1} can be equivalently reformulated as
\begin{align}\label{EQ:Problem2}
\min_{\M{x}}~F(\M{x})+R(\M{x}),~\mathrm{s.t.~} \sqrt{\nicefrac{1}{2}(\M{I}-\bar{\M{M}})}\M{x}=0.
\end{align}
To solve problem \eqref{EQ:Problem2}, we consider the saddle-point formulation
\begin{align}\label{EQ-Problem2-saddle}
\min_{\M{x}}\max_{\M{u}}\  \underbrace{F(\M{x})+R(\M{x})+\Big\langle \M{x},\sqrt{\nicefrac{1}{2}(\M{I}-\bar{\M{M}})}\M{u} \Big\rangle}_{:=L(\M{x},\M{u})},
\end{align}
where $\M{u}\in\mathbb{R}^{n\times d}$ is the dual variable. Denoting $\partial L(\M{x},\M{u})$ as the saddle-point subdifferential of $L(\M{x},\M{u})$, we have that
\begin{align}\label{l-mapping}
\partial L(\M{x},\M{u})=\left(
         \begin{array}{cc}
           \nabla F(\M{x})+\partial R(\M{x})+\sqrt{\nicefrac{1}{2}(\M{I}-\bar{\M{M}})}\M{u} \\
           -\sqrt{\nicefrac{1}{2}(\M{I}-\bar{\M{M}})}\M{x} \\
         \end{array}
       \right).
\end{align}
From \cite[Section 15.3]{Duality1}, it holds that $(\M{x}^\star,\M{u}^\star)$ solves problem \eqref{EQ-Problem2-saddle}, i.e., $\M{x}^\star$ solves problem \eqref{EQ:Problem2} and $\M{u}^\star$ solves its dual problem, iff $0\in \partial L(\M{x}^\star,\M{u}^\star)$. We split the smooth term $\nabla F(\M{x})$ and non-smooth term $\partial R(\M{x})$ in $\partial L(\M{x},\M{u})$ by introducing an auxiliary variable
$$
\M{z}:=\M{x}+\alpha \partial R(\M{x}).
$$
Based on this auxiliary variable, we define the mapping $\mathcal{M}(\M{x},\M{z},\M{u})$, which takes the following form
\begin{align}\label{m-mapping}
\mathcal{M}=\left(
                                     \begin{array}{c}
\M{z}-\M{x}+\alpha \nabla F(\M{x})+\alpha \sqrt{\nicefrac{1}{2}(\M{I}-\bar{\M{M}})} \M{u}\\
-\sqrt{\nicefrac{1}{2}(\M{I}-\bar{\M{M}})}\M{z}\\
(\M{x}+\alpha \partial R(\M{x}))-\left(\M{I}-{\nicefrac{1}{2}(\M{I}-\bar{\M{M}})}\right)\M{z}
                                     \end{array}
                                   \right).
\end{align}
Then, we have the following results.
\begin{prop}\label{Pro2}
Let $\M{z}=\M{x}+\alpha \partial R(\M{x})$. With Assumption \ref{MixingMatrix}, a point $(\M{x}^\star,\M{z}^\star,\M{u}^\star)$ solves $0\in\mathcal{M}(\M{x},\M{z},\M{u})$ if and only if $(\M{x}^\star,\M{u}^\star)$ solves $0\in \partial L(\M{x},\M{u})$.
\end{prop}
\begin{IEEEproof}
See Appendix \ref{AP:Pro2}.
\end{IEEEproof}
By Proposition \ref{Pro2}, we know that finding the zero point of $\partial L(\M{x},\M{u})$ is equivalent to finding the zero point of the mapping $\mathcal{M}(\M{x},\M{z},\M{u})$. From the structure of \eqref{m-mapping}, to find the zero point of $\mathcal{M}(\M{x},\M{z},\M{u})$, we develop the following iteration
\begin{align*}
\M{z}^t&=\M{x}^t-\alpha \nabla F(\M{x}^t)-\alpha\sqrt{\nicefrac{1}{2}(\M{I}-\bar{\M{M}})}\M{u}^t,\\
\M{u}^{t+1}&=\M{u}^{t}+\nicefrac{\theta_tp}{\alpha} \sqrt{\nicefrac{1}{2}(\M{I}-\bar{\M{M}})}\M{z}^t,\\
\M{x}^{t+1}&=\mathrm{prox}_{\alpha R}\left(\left(\M{I}-\nicefrac{\theta_t}{2}{(\M{I}-\bar{\M{M}})}\right)\M{z}^t\right).
\end{align*}
Setting scaled dual variable $\M{y}^t=\sqrt{\nicefrac{1}{2}(\M{I}-\bar{\M{M}})}\M{u}^t$, from the above iteration, we develop $\textsc{MG-Skip}$ as the following form
\begin{subequations}\label{MG-SKIP-NEW-ITER}
\begin{align}
\M{z}^t&=\M{x}^t-\alpha \nabla F(\M{x}^t)-\alpha\M{y}^t,\\
\bar{\M{z}}^{t+1}&=\nicefrac{\theta_t}{2}(\M{I}-\bar{\M{M}})\M{z}^t,\\
\M{y}^{t+1}&=\M{y}^{t}+\nicefrac{p}{\alpha} \bar{\M{z}}^{t+1},\label{MG-SKIP-NEW-ITER-dual}\\
\M{x}^{t+1}&=\mathrm{prox}_{\alpha R}\left(\M{z}^t-\bar{\M{z}}^{t+1}\right),
\end{align}
\end{subequations}
where $\alpha>0$ is the stepsize and $\theta_t=1$ with probability $p$ and $\theta_t=0$ with probability $1-p$. Splitting the updates to the nodes, we elaborate the decentralized implementation of $\textsc{MG-Skip}$ \eqref{MG-SKIP-NEW-ITER} in Algorithm \ref{DP-ProxPro}. To reduce communication costs, similar as ProxSkip \cite{ProxSkip}, we trigger the communication with a probability $p$, meaning that not every iteration necessitates communication. In expectation the communication is occurred every $\nicefrac{1}{p}$ iterations, which can be very rare if $p$ is small. Additionally, we employ multi-gossip communication, which mitigates the limitations of probabilistic local updates by not relying on specific network connectivity.

\begin{algorithm}[!t]
\caption{$\textsc{MG-Skip}$}
\label{DP-ProxPro}
    \SetAlgoNlRelativeSize{0}
     \textbf{Input:} $\alpha\in(0,\nicefrac{2}{L})$ and $p\in(0,1]$. $\M{W}=[W_{ij}]_{i,j=1}^n$\\
     \textbf{Initial:} $\m{x}^0_{i}\in\mathbb{R}^d$ and $\m{y}^0_i = 0\in\mathbb{R}^d,i=1,\ldots,n$\\
     Flip coins $[\theta_0,\ldots,\theta_{T-1}] \in [0,1]^{T}$ with $\mathop{\M{P}}(\theta_t =1) = p$\\
    \For{$t=0,1,\dotsc,T-1$, every node $i$}{
    $\m{z}^t_i=\m{x}^t_i-\alpha\nabla f_i(\m{x}^t_i)-\alpha \m{y}^t_i$ \hfill \textcolor[rgb]{0.3529,0.3961,0.4235}{$\triangleright$ update the variate $\M{z}^t_i$}\\
      \eIf{$\theta_t=1$}{
                                    $\bar{\m{z}}^{k+1}_i=\nicefrac{1}{2}\textsc{FastGoss}([W_{ij}]_{j\in\mathcal{N}_i},\{{\m{z}}_j^t\}_{j\in\mathcal{N}_i},K)$\\
                                    $\m{y}^{t+1}_i=\m{y}^{t}_i+\nicefrac{p}{\alpha}\bar{\m{z}}^{k+1}_i$\\
                                    $\m{x}^{t+1}_i=\mathrm{prox}_{\alpha r}\left({\m{z}}^t_i-\bar{\m{z}}^{k+1}_i\right)$\\
      }{
      $\m{x}^{t+1}_i=\mathrm{prox}_{\alpha r}({\m{z}}^{t}_i),\m{y}^{t+1}_i =  \m{y}^{t}_i$ \hfill \textcolor[rgb]{0.3529,0.3961,0.4235}{$\triangleright$ skip Comm.}
      }
    }
\textbf{Output:} $\M{x}_i^T,~i=1,\dots,n$\\
\vspace{2mm}
\textbf{procedure} $\textsc{FastGoss}([W_{ij}]_{j\in\mathcal{N}_i},\{\M{s}_j\}_{j\in\mathcal{N}_i},K)$\\
$\eta=\nicefrac{\big(1-\sqrt{1-\rho^2}\big)}{\big(1+\sqrt{1+\rho^2}\big)}$, $\m{s}_i^{0}=\m{s}_i^{-1}={\m{z}}^{t}_i$\\
\For{$k=1,\dotsc,K-1$, every node $i$}{
$\m{s}_i^{k+1}=(1+\eta)\sum_{j\in \mathcal{N}_i}W_{ij}\m{s}_j^{k}-\eta\m{s}_i^{k-1}$
}
\textbf{return} $(\m{s}_i^0-\m{s}_i^K),i=1,\ldots,n$
\end{algorithm}

\subsection{Relation with Existing Algorithmic Frameworks}
Let $\Delta_F^t:= \nabla F(\M{x}^{t-1})-\nabla F(\M{x}^t)$. Eliminating the scaled dual variable $\M{y}^t$ on \eqref{MG-SKIP-NEW-ITER}, $\textsc{MG-Skip}$ takes the following updates.
\begin{align*}
\M{z}^t&=\left(\M{I}-\nicefrac{p\theta_{t-1}}{2}(\M{I}-\bar{\M{M}})\right){\M{z}}^{t-1}+\left(\M{x}^t-\M{x}^{t-1}\right)+\alpha \Delta_F^t,\\
\M{x}^{t+1}&=\m{prox}_{\alpha R}\left(\left(\M{I}- \nicefrac{\theta_t}{2}(\M{I}-\bar{\M{M}})\right){\M{z}}^{t}\right).
\end{align*}
Recall the framework PUDA \cite{Sulaiman2021}, whose primal updates take the following form by eliminating its dual variable
\begin{subequations}\label{PUDA-UP}
\begin{align}
{\M{z}}^t&=\mathcal{B}{\M{z}}^{t-1}+\mathcal{C}(\M{x}^{t}-\M{x}^{t-1})
+\alpha\Delta_F^t,\\
\mathbf{x}^{t+1}&=\mathrm{prox}_{\alpha R}(\mathcal{A}{\M{z}}^t),
\end{align}
\end{subequations}
where $\mathcal{A,B,C}$ are suitably chosen matrices. The convergence condition of PUDA is $\mathcal{A}^2\preceq\mathcal{B}\prec \M{I}$ and $0\preceq\mathcal{C}\preceq 2\M{I}$. Then, we recall the primal updates of $\mathcal{ABC}$-framework \cite{Xu2021}
\begin{subequations}\label{ABC-UP}
\begin{align}
\M{w}^{t}&=\mathcal{B}\M{w}^{t-1}+\mathcal{C}(\M{x}^{t}-\M{x}^{t-1})+\alpha\mathcal{A}\Delta_F^t,\\
\mathbf{x}^{t+1}&=\mathrm{prox}_{\alpha R}(\mathbf{w}^{t}).
\end{align}
\end{subequations}
Letting $\M{w}^{t}:=\mathcal{A}{\M{z}}^{t}$ in PUDA \eqref{PUDA-UP}, it holds that PUDA \eqref{PUDA-UP} and $\mathcal{ABC}$-framework \eqref{ABC-UP} are equivalent, since they generate the same sequence $\{\M{x}^t\}$. We first consider the special case of $\textsc{MG-Skip}$ with $p=1$. In this case, choosing $\mathcal{A}=\mathcal{B}=\M{I}-\nicefrac{1}{2}(\M{I}-\bar{\M{M}})$ and $\mathcal{C}=\M{I}$ in \eqref{PUDA-UP} and $\mathcal{A}=\mathcal{B}=\mathcal{C}=\M{I}-\nicefrac{1}{2}(\M{I}-\bar{\M{M}})$ in \eqref{ABC-UP}, $\textsc{MG-Skip}$ can be derived from PUDA \cite{Sulaiman2021} and $\mathcal{ABC}$-framework \cite{Xu2021}, respectively.

When $p\in(0,1)$, by comparison, since $\mathcal{A,B,C}$ are required to be time-invariant and deterministic in PUDA and $\mathcal{ABC}$-framework, $\textsc{MG-Skip}$ can not be unified with them. In addition, if we only consider the one-step iteration, by choosing $\mathcal{A}=\M{I}- \nicefrac{\theta_t}{2}(\M{I}-\bar{\M{M}})$, $\mathcal{B}= \M{I}-\nicefrac{p\theta_{t-1}}{2}(\M{I}-\bar{\M{M}})$, $\mathcal{C}=\M{I}$ in \eqref{PUDA-UP}, the updates of $\textsc{MG-Skip}$ at $t$-th iteration can be derived. However, despite this, the one-step iteration of $\textsc{MG-Skip}$ does not satisfy the convergence condition of PUDA. Specifically, when $\theta_{t-1} = 1$ and $\theta_t = 0$, we have $\mathcal{A} = \M{I}$ and $\mathcal{B} = \M{I} - \nicefrac{p}{2}(\M{I}-\bar{\M{M}})$, i.e., $\mathcal{A}^2\succ \mathcal{B}$. Therefore, the existing convergence analyses in \cite{Sulaiman2021} and \cite{Xu2021} are not applicable to $\textsc{MG-Skip}$.

\section{Main results}\label{SEC4}
In this section, we provide the convergence analysis for $\textsc{MG-Skip}$ (Algorithm \ref{DP-ProxPro}). Based on this convergence result, we further discuss the relation and improvement of the proposed $\textsc{MG-Skip}$ with existing algorithms.

Before presenting our results, we give the optimality condition for problems \eqref{EQ:Problem1} and \eqref{EQ:Problem2}.
\begin{lem}\label{LEM-1}
Suppose that Assumptions \ref{ASS1} and \ref{MixingMatrix} hold. If there exists a point $(\M{x}^{\star},\M{z}^{\star},\M{u}^{\star})$ that satisfies:
\begin{subequations}\label{KKT-Condition}
\begin{align}
\M{z}^{\star}&=\M{x}^{\star}-\alpha \nabla F(\M{x}^{\star})-\alpha \sqrt{\nicefrac{1}{2}(\M{I}-\bar{\M{M}})} \M{u}^{\star},\label{KKT-Condition-1}\\
0&=\sqrt{\nicefrac{1}{2}(\M{I}-\bar{\M{M}})}\M{z}^{\star},\label{KKT-Condition-2}\\
\M{x}^{\star}&=\mathrm{prox}_{\alpha R}(\M{z}^{\star}), \label{KKT-Condition-3}
\end{align}
\end{subequations}
it holds that $\M{x}^{\star}=[\m{x}^{\star},\m{x}^{\star},\cdots,\m{x}^{\star}]\tr$, where $\m{x}^{\star}\in \mathbb{R}^d$ is a solution to problem \eqref{EQ:Problem1}.
\end{lem}
\begin{IEEEproof}
See Appendix \ref{AP:LEM-1}.
\end{IEEEproof}
Under Assumptions \ref{ASS1} and \ref{MixingMatrix}, from \cite[Lemma 1]{Sulaiman2021} or \cite[Lemma 5]{Xu2021}, we can directly derive that a point $(\M{x}^{\star},\M{z}^{\star},\M{u}^{\star})$ satisfying \eqref{KKT-Condition} exists and is unique. Moreover, it is not difficult to verify that, when $\theta_t=1$, the fixed point of $\textsc{MG-Skip}$ \eqref{MG-SKIP-NEW-ITER}, i.e., $(\M{x}^{t+1},\M{u}^{t+1})=(\M{x}^{t},\M{u}^{t})$, satisfies the optimality condition \eqref{KKT-Condition}. It implies that if the sequence $\{\M{x}^t\}_{t\geq0}$ generated by $\textsc{MG-Skip}$ \eqref{MG-SKIP-NEW-ITER} is convergent, it will converges to $\M{x}^{\star}=[\m{x}^{\star},\m{x}^{\star},\cdots,\m{x}^{\star}]\tr$, where $\m{x}^{\star}\in \mathbb{R}^d$ is a solution to problem \eqref{EQ:Problem1}.

Subsequently, we present the following lemma, demonstrating that after a single step of $\textsc{MG-Skip}$ \eqref{MG-SKIP-NEW-ITER}, the distance between $(\M{x}^{t+1},\M{u}^{t+1})$ and $(\M{x}^\star,\M{u}^\star)$ can be bounded in terms of the distances $\|(\M{x}^t-\alpha\nabla F(\M{x}^t))-(\M{x}^\star-\alpha\nabla F(\M{x}^\star))\|^2$ and $\|\M{u}^t-\M{u}^\star\|^2$. This nonexpansiveness plays a crucial role in the ensuing proof.

\begin{lem}\label{LEM-2}
Suppose that Assumptions \ref{ASS1} and \ref{MixingMatrix} hold, $0<p\leq1$ and $\mu\geq0$. Let $\{(\M{x}^t,\M{z}^t,\M{u}^t)\}_{t\geq0}$ be the sequence generated by $\textsc{MG-Skip}$ \eqref{MG-SKIP-NEW-ITER}, where $\M{u}^t$ satisfying that $\M{y}^t=\sqrt{\nicefrac{1}{2}(\M{I}-\bar{\M{M}})}\M{u}^t$ is the dual variable.
For any $(\M{x}^{\star},\M{z}^{\star},\M{u}^{\star})$ satisfying $0\in\mathcal{M}(\M{x}^{\star},\M{z}^{\star},\M{u}^{\star})$,
it holds that
\begin{align}\label{PROOF-IEQ8}
&\mathbb{E}\left[\left\|\M{x}^{t+1}-\M{x}^\star \right\|_{\m{F}}^2~|~\theta_t\right]+\frac{\alpha^2}{p^2}\mathbb{E}\left[\left\|\M{u}^{t+1}-\M{u}^\star \right\|_{\m{F}}^2~|~\theta_t\right]\nonumber\\
&\leq  \Big\|\M{v}^t-\M{v}^{\star} \Big\|_{\m{F}}^2+\Big(\frac{\alpha^2}{ p^2}-\frac{\alpha^2}{5}\Big) \Big\|\M{u}^t-\M{u}^{\star} \Big\|_{\m{F}}^2,
\end{align}
where $\M{v}^t=\M{x}^t-\alpha \nabla F(\M{x}^t)$ and $\M{v}^{\star}=\M{x}^{\star}-\alpha \nabla F(\M{x}^{\star})$.
\end{lem}
\begin{IEEEproof}
See Appendix \ref{AP:LEM-2}.
\end{IEEEproof}
\subsection{Convergence Analysis}
\begin{thm}\label{THM-2}
Suppose that Assumptions \ref{ASS1} and \ref{MixingMatrix} hold, $0<p\leq1$ and $\mu>0$. Let $\{(\M{x}^t,\M{z}^t,\M{u}^t)\}_{t\geq0}$ be the sequence generated by $\textsc{MG-Skip}$ \eqref{MG-SKIP-NEW-ITER}, where $\M{u}^t$ satisfying that $\M{y}^t=\sqrt{\nicefrac{1}{2}(\M{I}-\bar{\M{M}})}\M{u}^t$ is the dual variable.
For any $(\M{x}^{\star},\M{z}^{\star},\M{u}^{\star})$ satisfies $0\in\mathcal{M}(\M{x}^{\star},\M{z}^{\star},\M{u}^{\star})$, define the Lyapunov function,
$$
\Psi^{t}=\left\|\M{x}^t-\M{x}^{\star}\right\|_{\m{F}}^2+\nicefrac{\alpha^2}{ p^2}\left\|\M{u}^t-\M{u}^{\star}\right\|_{\m{F}}^2.
$$
If $0<\alpha<\frac{2}{L}$,
it holds that
\begin{align}\label{LinearRate1}
\mathbb{E}\left[\left\|\M{x}^t-\M{x}^{\star}\right\|_{\m{F}}^2\right]\leq \zeta^{t}\Psi^{0},
\end{align}
where
$$
{\zeta}:=\max\left\{(1-\alpha\mu)^2,(1-\alpha L)^2,1-\nicefrac{ p^2}{5}\right\}<1.
$$
Moreover, it holds that
\begin{align*}
\lim_{t\rightarrow\infty}\left\|\M{x}^t-\M{x}^\star\right\|_{\m{F}}^2&=0, \text{ almost surely}, \\
\lim_{t\rightarrow\infty}\left\|\M{u}^t-\M{u}^\star\right\|_{\m{F}}^2&=0, \text{ almost surely}.
\end{align*}
\end{thm}
\begin{IEEEproof}
See Appendix \ref{AP:THM-2}.
\end{IEEEproof}

Theorem \ref{THM-2} shows the linear convergence of $\textsc{MG-Skip}$. Since we use multi-gossip communications when the communication is triggered, which is equivalent to an information fusion on a sufficiently well-connected network, the linear convergence rate $\zeta$ remains unaffected by the network topology, which is consistent with \cite{Xu2021,Scaman2017,Yuan2021}.

In addition, leveraging Theorem \ref{THM-2} as a foundation, we proceed to offer a more in-depth analysis.

\noindent\textbf{$\bullet$ $\textsc{MG-Skip}$ allows for skipping of partial communication for free without assuming that the network is sufficiently well-connected.}
Let
$$\zeta_0:=\max\{(1-\alpha\mu)^2,(1-\alpha L)^2\}.$$
Since $\zeta_0\leq 1-\mu\alpha$, setting $\alpha=\nicefrac{1}{5L}$, from \eqref{LinearRate1}, we have
$$\mathbb{E}\left[\left\|\M{x}^t-\M{x}^{\star}\right\|_{\m{F}}^2\right]\leq \left(\max\{1-\nicefrac{1}{5\kappa},1-\nicefrac{p^2}{5}\}\right)^{t}\Psi^{0},$$
i.e., the convergence rate of $\textsc{MG-Skip}$ is bounded by
$$
\mathrm{RATE}_{\textsc{MG-Skip}}=\mathcal{O}\left(\max\left\{1-\nicefrac{1}{\kappa},1-p^2\right\}\right).
$$
Therefore, it holds that
$$
\mathrm{RATE}_{\textsc{MG-Skip}}\equiv \mathcal{O}\left(1-\nicefrac{1}{\kappa}\right),\ \text{when } p\in[\nicefrac{1}{\sqrt{\kappa}},1].
$$
Thus, it implies that, when setting $p\in[\nicefrac{1}{\sqrt{\kappa}},1]$, it is beneficial for saving communication costs. More importantly, note that the convergence rate remains unchanged as we decrease $p$ from $1$ to $\nicefrac{1}{\sqrt{\kappa}}$. This is the reason why we can often skip the $\textsc{FastGoss}(\cdot)$, and get away with it for free, i.e., without any deterioration of the convergence rate. Recall the convergence rate of ProxSkip \cite{ProxSkip} when $\alpha=\mathcal{O}(\nicefrac{1}{L})$, which is
$$
\mathrm{RATE}_{\text{ProxSkip}}=\mathcal{O}\left(\max\{1-\nicefrac{1}{\kappa},1-p^2(1-\rho)\}\right).
$$
When assuming the network is sufficiently well-connected, i.e., $\rho\leq1-\nicefrac{1}{\kappa}$, it holds that
$$
\mathrm{RATE}_{\text{ProxSkip}}\equiv \mathcal{O}\left(1-\nicefrac{1}{\kappa}\right),\ \text{when } p\in[\nicefrac{1}{\sqrt{\kappa(1-\rho)}},1],
$$
which implies that ProxSkip \cite{ProxSkip} is free to skip communications only when the network is sufficiently well-connected.\\

\noindent \textbf{$\bullet$ $\textsc{MG-Skip}$ attains the optimal communication complexity for composite problem \eqref{EQ:Problem2}.} Noting that $\zeta_0\leq(1-\alpha\mu)$, from \eqref{LinearRate1}, we have
\begin{align*}
\mathbb{E}\left[\left\|\M{x}^t-\M{x}^{\star}\right\|_{\m{F}}^2\right]\leq \left(\max\{1-\alpha\mu,1-\nicefrac{ p^2}{5}\}\right)^t\Psi^{0}.
\end{align*}
Setting $\alpha=\mathcal{O}(\nicefrac{1}{L})$, we have
$$
\max\{1-\alpha\mu,1-\nicefrac{p^2}{5}\}=\max\{\mathcal{O}(1-\nicefrac{1}{\kappa}),\mathcal{O}(1-p^2)\}.
$$
Thus, it holds that the convergence rate $\zeta$ of $\textsc{MG-Skip}$ can be bounded as follows
$$
\zeta\leq\max\{\mathcal{O}(1-\nicefrac{1}{\kappa}),\mathcal{O}(1-p^2)\}.
$$
For any given target accuracy $\epsilon>0$, by solving
$$
\Exp{\|\M{x}^t-\M{x}^{\star}\|^2_{\m{F}}}\leq\epsilon \Psi^{0},
$$
the computation complexity $T$ (the number of gradient computations) is
$$
T= \max\left\{\mathcal{O}\big(\kappa\log\nicefrac{1}{\epsilon}\big),\mathcal{O}\big(\nicefrac{1}{p^2}\log\nicefrac{1}{\epsilon}\big)\right\}.
$$
Since in each iteration we evaluate the $\textsc{FastGoss}(\cdot)$ with probability $p$, which requires $K=\lfloor\nicefrac{1}{\sqrt{1-\rho}}\rfloor$ consecutive rounds of communications, the expected communication complexity $C$ (the number of gradient communications) is
\begin{align*}
C &= p\left\lfloor\frac{1}{\sqrt{1-\rho}}\right\rfloor T\\
&\thickapprox\max\left\{\mathcal{O}\left(\frac{p\kappa}{\sqrt{1-\rho}}\log\nicefrac{1}{\epsilon}\right),\mathcal{O}\left(\frac{1}{p\sqrt{1-\rho}}\log\nicefrac{1}{\epsilon}\right)\right\}.
\end{align*}
Thus, when reaching $\Exp{\|\M{x}^t-\M{x}^{\star}\|^2_{\m{F}}}\leq\epsilon \Psi^{0}$, by choosing $p\in[\nicefrac{1}{\sqrt{\kappa}},1]$, the computation complexity $T$ and the communication complexity $C$ are
\begin{align*}
&T=\mathcal{O}\left(\kappa\log \nicefrac{1}{\epsilon}\right),\quad C=\mathcal{O}\left(\frac{p\kappa}{\sqrt{1-\rho}}\log \nicefrac{1}{\epsilon}\right).
\end{align*}
When $p\in[\nicefrac{1}{\sqrt{\kappa}},1]$, the computation complexity $T$ is independent of $p$, so decreasing the probability of communication (skipping part of the communication) does not increase the number of gradient computations required to achieve the target accuracy. On the other hand, since $p\leq 1$, $\textsc{MG-Skip}$ can achieve provable communication acceleration. Moreover, setting $p=\nicefrac{1}{\sqrt{\kappa}}$, it holds that
\begin{align*}
T=\mathcal{O}\left(\kappa\log\nicefrac{1}{\epsilon}\right),\quad
C=\mathcal{O}\left(\sqrt{\frac{\kappa}{1-\rho}}\log\nicefrac{1}{\epsilon}\right).
\end{align*}
From \cite{Scaman2017,Kovalev2020,HuanLi2020,HuanLi2022}, $\textsc{MG-Skip}$ achieves the optimal communication complexity. To obtain this complexity, the global parameters $\kappa$ and $\rho$ are used. To calculate $\kappa$,
we take variables $(r^1_i,s^1_i)=(L_i,\mu_i),i\in\mathcal{V}$ and update them by
$$
(r_i^{k+1},s_i^{k+1})=(\max\{r_i^k,r_j^k:j\in\mathcal{N}_i\},\min\{s_i^k,s_j^k:j\in\mathcal{N}_i\}).
$$
The global parameters $L=r_i^{n}$, $\mu=s_i^{n}$, and $\kappa=\nicefrac{r_i^{n}}{s_i^{n}}$ can be calculated by agent $i$ within $n-1$ steps. When the network topology and weights are known in advance, as in \cite{Xu2021,Scaman2017,Yuan2021}, we can assume that the global parameter $\rho$ is known and accessible by agents. Conversely, one can estimate $\rho$ by \cite[Proposition 5]{Nedic2018}, where estimates of $\rho$ for many types of networks are given.\\

\noindent \textbf{$\bullet$ $\textsc{MG-Skip}$ is the first algorithm with both random communication skippings and linear convergence for decentralized composite optimization problem \eqref{EQ:Problem1}.} In \cite{RandProx}, a primal-dual unified framework with randomized dual update, named RandProx, is given for minimizing $f(x)+g(x)+h(Ax)$, where $f$, $g$ and $h$ are proper closed convex, $f$ is smooth, and $A$ is a bounded linear operator. In terms of \cite[Example 2]{RandProx}, the probability local updates mechanism can be viewed as a special case of its proposed randomized dual update. Thus, RandProx can be used to solve the decentralized composite optimization problem \eqref{EQ:Problem1} with local updates. From \cite[Eq. (28)]{RandProx}, the sequence $\{(\M{x}^{t},\M{u}^t)\}_{\{t\geq0\}}$ generated by RandProx satisfies that
\begin{align*}
&\mathbb{E}\left[\left\|\M{x}^{t+1}-\M{x}^\star \right\|_{\m{F}}^2\right]+\nicefrac{\alpha^2}{p^2}\mathbb{E}\left[\left\|\M{u}^{t+1}-\M{u}^\star \right\|_{\m{F}}^2\right]\nonumber\\
&\leq  \zeta_0\left\|\M{x}^t-\M{x}^{\star}\right\|_{\m{F}}^2+ \nicefrac{\alpha^2}{ p^2}\left\|\M{u}^t-\M{u}^{\star} \right\|_{\m{F}}^2,
\end{align*}
which fails to deduce a linear convergence rate. In fact, from the unified proof of RandProx for general composite optimization $\min_{x} \{f(x)+g(x)+h(Ax)\}$, i.e., \cite[Appendix D]{RandProx}, we have
\begin{align*}
X^{t+1}+U^{t+1}\leq r_1 X^{t} + r_2U^t,
\end{align*}
where $X^{t}$ and $U^t$ measure the distance of the primal and dual states generated by RandProx at $t$-th iteration to the solution of the primal and dual problems, respectively,
$r_1:=\frac{\zeta_0}{1+\gamma \mu_g}$, and
$r_2:=(\frac{1+\omega}{\tau}+2 \omega L_{h}) / (\frac{1+\omega}{\tau}+2(1+\omega) L_{h})$,
where the $\gamma$ and $\tau$ are the primal and dual step sizes, $\omega$ is the relative variance of the dual update, and $L_{h}$ is the smooth factor of $h$. Without assuming that $h$ is smooth, we have $r_2=1$, i.e.,
\begin{align*}
X^{t+1}+U^{t+1}\leq r_1 X^{t} + U^t \nrightarrow \text{linear convergence}.
\end{align*}
Therefore, when $g\neq0$, the assumption that ``$h$ is smooth'' is important to the proof presented in \cite{RandProx}. However, to the considered problem \eqref{EQ:Problem2}, which can be equivalently written as
$\min_{x}\{ f(x)+g(x)+h(Ax)\}$, where $f(x)=F(\M{x})$, $g(x)=R(\M{x})\neq0$, and $h(Ax)=\iota_{\{0\}}(\sqrt{\nicefrac{1}{2}(\M{I}-\bar{\M{M}})}\M{x})$ is non-smooth, with $\iota_{\{0\}}$ being the indicator function of $\{0\}$, with $\iota_{\{0\}}:x\mapsto 0, \text{ if }x=0;\ x\mapsto \infty, \text{otherwise}$. Thus, the linear convergence results established in \cite{RandProx} are not applicable to it. To $\textsc{MG-Skip}$, by introducing the auxiliary variable $\M{z}=\M{x}+\alpha \partial R(\M{x})$, we split the smooth term and non-smooth of the primal updates as
\begin{align*}
\M{z}^t&=\M{x}^t-\alpha \nabla F(\M{x}^t)-\alpha\sqrt{\nicefrac{1}{2}(\M{I}-\bar{\M{M}})}\M{u}^t,\\
(\M{I}+\alpha\partial R)(\M{x}^{t+1})&=\left(\M{I}- {\nicefrac{\theta_t}{2}(\M{I}-\bar{\M{M}})}\right)\M{z}^t.
\end{align*}
Based on the separable structure of $\textsc{MG-Skip}$, we first establish an upper bound for the smooth part $\mathbb{E}[\|(\M{I} - {\nicefrac{\theta_t}{2}(\M{I}-\bar{\M{M}})}) \M{z}^t - \M{z}^\star\|^2_{\m{F}}]$, which includes the term $-\nicefrac{\alpha^2}{5}\left\|\M{u}^t - \M{u}^\star\right\|_{\m{F}}^2$. Then, utilizing the nonexpansivity of the proximal mapping, we demonstrate that
\begin{align*}
&\mathbb{E}\left[\left\|\M{x}^{t+1} - \M{x}^\star\right\|_{\m{F}}^2\right] + \nicefrac{\alpha^2}{p^2} \mathbb{E}\left[\left\|\M{u}^{t+1} - \M{u}^\star\right\|_{\m{F}}^2\right] \\
&\leq \zeta_0 \left\|\M{x}^t - \M{x}^\star\right\|_{\m{F}}^2 + (1 - \nicefrac{p^2}{5}) \left(\nicefrac{\alpha^2}{p^2} \left\|\M{u}^t - \M{u}^\star\right\|_{\m{F}}^2\right),
\end{align*}
which indicates a linear convergence rate.

\subsection{Further Discussions}
\begin{figure*}[!t]
  \centering
  \setlength{\abovecaptionskip}{-2pt}
  \includegraphics[width=0.95\linewidth]{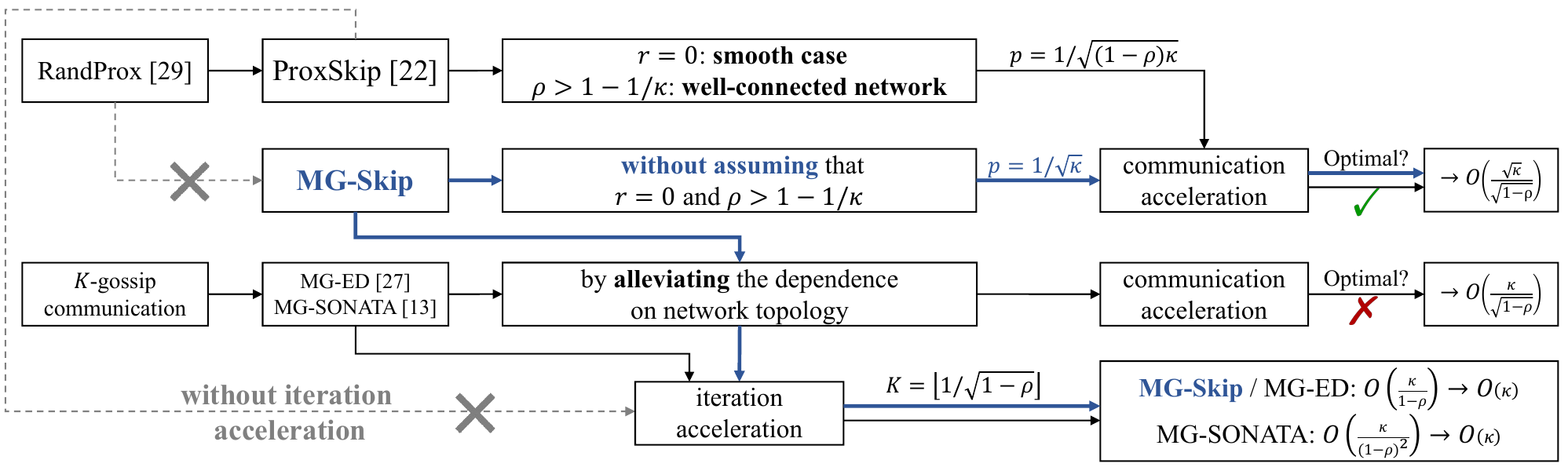}
  \caption{The relation and improvements of $\textsc{MG-Skip}$ compared to exiting decentralized algorithms.}
  \label{imp}
\end{figure*}

\emph{Local updates}, including periodic local updates \cite{Liu2023,Alghunaim2023} and probabilistic local updates \cite{ProxSkip}, are recently emerging mechanisms for accelerating communication.
However, their efficacy lies in accelerating communication complexity solely under the conditions of a smooth loss function and a sufficiently well-connected network ($\rho\leq 1-\nicefrac{1}{\kappa}$). In addition, \cite{Liu2023,Alghunaim2023,ProxSkip} are designed for the smooth case, thus they cannot be used to solve the composite decentralized optimization we consider. Take ProxSkip for example
\begin{align*}
\tilde{\mathcal{O}}\left(\frac{\kappa}{1-\rho}\right) \longrightarrow\tilde{\mathcal{O}}\left(\frac{\sqrt{\kappa}}{\sqrt{1-\rho}}\right), \text{ for } r=0 \text{ and }\rho\leq 1-\nicefrac{1}{\kappa}.
\end{align*}

\emph{$K$-gossip} (multi-gossip) is a common acceleration mechanism in decentralized optimization. In \cite{Scaman2017,HuanLi2020,HuanLi2022}, combining with Nesterov's accelerated gradient descent methods \cite{Nesterov}, it has been used to develop optimal algorithms for smooth and strongly convex decentralized optimization in networks. In \cite{Yuan2021}, $K$-gossip has been used to alleviate network topology dependence of D-SGD \cite{Stich2020}. For composite settings, $K$-gossip has been adopted to accelerate algorithm in \cite{SONATA,PMGT-VR}. However, $K$-gossip exclusively serves to remove network topology dependence and falls short of attaining acceleration at the problem level. In other words, the complexity's reliance on $\kappa$ remains unchanged irrespective of the deployment of this mechanism. Take SONOTA and PMGT-VR for example
\begin{align*}
\left\{\begin{array}{l}
\text{\# of Comm.}: \tilde{\mathcal{O}}\left(\frac{1}{(1-\rho)^2}\right) \xrightarrow{K\text{-gossip}}\tilde{\mathcal{O}}\left(\frac{1}{\sqrt{1-\rho}}\right),\\
\text{\# of Comm.}: \tilde{\mathcal{O}}(\kappa) \xrightarrow{K\text{-gossip}} \tilde{\mathcal{O}}(\kappa), \text{ \textbf{unchanged}}.
\end{array}\right.
\end{align*}

In this work, we combine $K$-gossip with local updates and propose a new decentralized implementation mechanism $\textsc{MG-Skip}$. Surprisingly, combining $K$-gossip and local updates compensates for the shortcomings of both sides and achieves a provable double acceleration in communication complexity, i.e.,
\begin{align*}
\left\{\begin{array}{l}
\text{\# of Comm.}: \tilde{\mathcal{O}}\left(\frac{1}{1-\rho}\right) \longrightarrow\tilde{\mathcal{O}}\left(\frac{1}{\sqrt{1-\rho}}\right),\\
\text{\# of Comm.}: \tilde{\mathcal{O}}(\kappa) \longrightarrow \tilde{\mathcal{O}}(\sqrt{\kappa}), \text{ \textbf{improved}},
\end{array}\right.
\end{align*}
without assuming $r=0$ and $\rho\leq 1-\nicefrac{1}{\kappa}$. We show the relation and improvements of $\textsc{MG-Skip}$ compared to exiting algorithms as Fig. \ref{imp}. In summary, $\textsc{MG-Skip}$ circumvents the limitations of ProxSkip \cite{ProxSkip} by not requiring specific network connectivity. Instead, it achieves the ``equivalent" connectivity by performing a multi-gossip communication when the communication is triggered. Furthermore, the provable linear convergence of $\textsc{MG-Skip}$ affirms the benefits of local updates in the decentralized composite setup.

\begin{figure*}[!t]
  \centering
  \setlength{\abovecaptionskip}{-2pt}
  \subfigure{
  \includegraphics[width=0.48\linewidth]{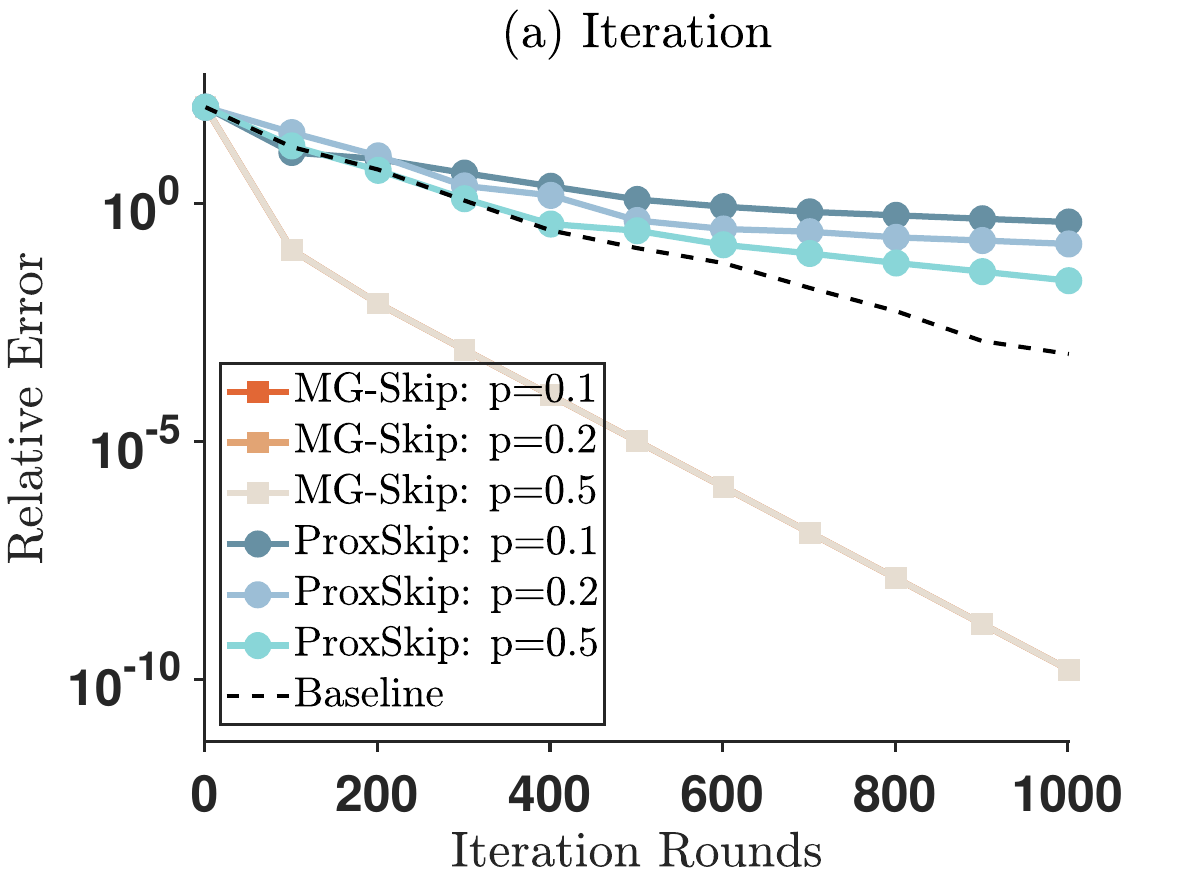}}
  \subfigure{
  \includegraphics[width=0.48\linewidth]{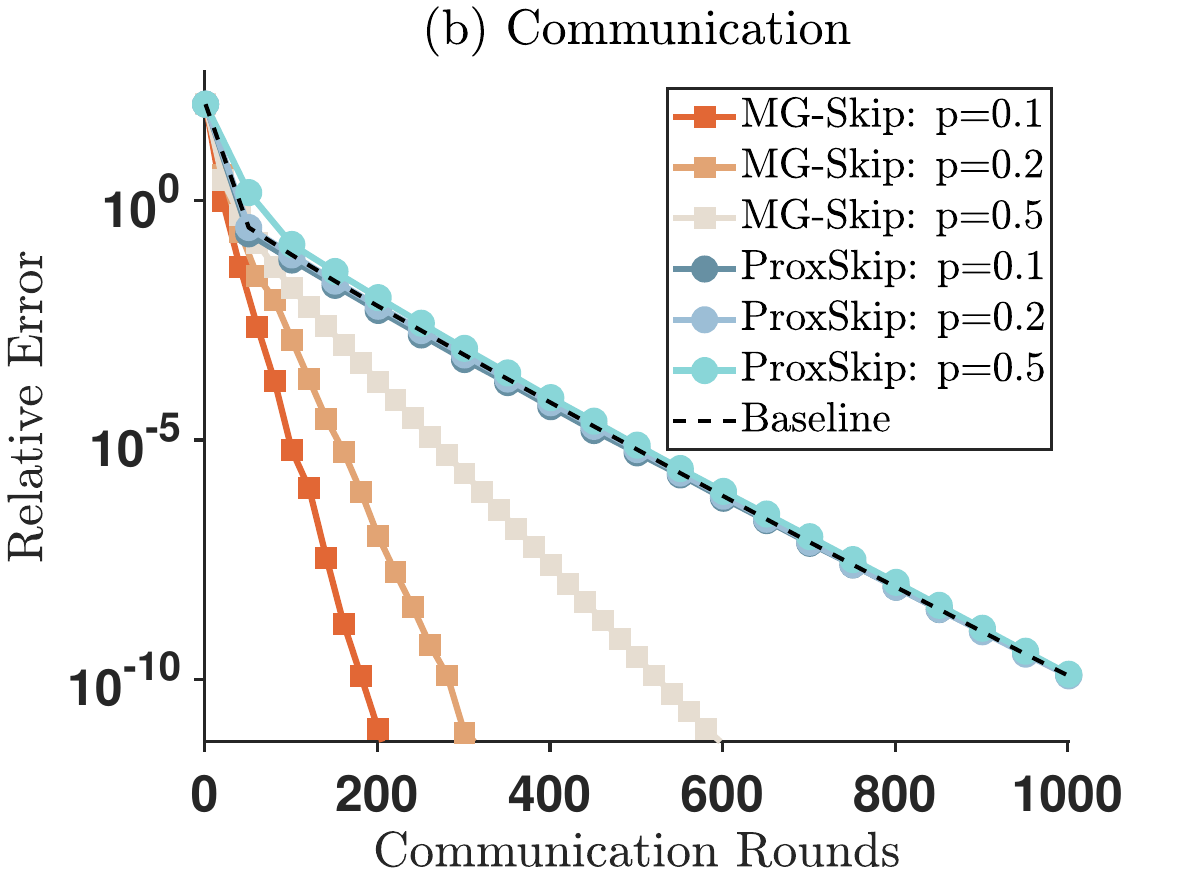}}
  \caption{SYNTHETIC convex function over $15$ nodes with $\rho\approx0.9424$. Here, baseline denotes ProxSkip $(p=1)$, which is equivalent to NIDS \cite{NIDS}/ED \cite{ExactDiffusion}).}
  \label{CASE1}
\end{figure*}

\section{Experimental Results}\label{SEC5}
In this section, we conduct several numerical experiments to validate the obtained theoretical results. All the algorithms are implemented in Matlab R2020b in a computer with 3.30 GHz AMD Ryzen 9 5900HS with Radeon Graphics and 16 GB memory.

\subsection{Decentralized Least Squares Problem on Synthetic Dataset}
In \cite{Alghunaim2023} and \cite{ProxSkip}, it is demonstrated that algorithms like those mentioned in \cite{Liu2023,Alghunaim2023,ProxSkip}, which incorporate local updates, can achieve communication acceleration exclusively when the network exhibits adequate connectivity, specifically when $\rho\leq 1-\nicefrac{1}{\kappa}$. In this experiment, our focus is on the scenario where $\rho\geq 1-\nicefrac{1}{\kappa}$, aiming to demonstrate that $\textsc{MG-Skip}$ can also achieve communication acceleration even when the network lacks sufficient connectivity.

We construct the decentralized least squares objective with
$f_i(\m{x})=\frac{1}{2}\|\m{A}_i\m{x}-\m{b}_i\|^2$ with $\mu \M{I}_d\preceq\m{A}_i\tr\m{A}_i\preceq L\M{I}_d$ and sample each $\m{b}_i\thicksim\mathcal{N}(0,\M{I}_d)$ for each node $i\in[n]$. We use a ring topology with 15 nodes for this experiment, where $\rho\approx0.9424$, and set $\kappa=\nicefrac{L}{\mu}=\nicefrac{0.5}{(1-\rho)}$, i.e., $\rho=1-\nicefrac{0.5}{\kappa}\geq 1-\nicefrac{1}{\kappa}$. Here, our performance comparison is solely between $\textsc{MG-Skip}$ and ProxSkip \cite{ProxSkip}. This approach is taken because, as stated in \cite[Section 3]{Alghunaim2023}, LED \cite{Alghunaim2023} and ProxSkip \cite{ProxSkip} are considered equivalent from an expectation perspective. The primary distinction lies in the selection of stepsize and the mixing matrix $\M{W}$. To all algorithms we use the same stepsies $\alpha=\nicefrac{1}{5L}$. The comparison results are shown in Fig. \ref{CASE1}. The relative error $\nicefrac{\|\M{x}^t-\M{x}^{\star}\|}{\|\M{x}^\star\|}$ is shown on the $y$-axis. Based on the results presented in Fig. \ref{CASE1}, it becomes evident that, in the case of a sparse network, such as a ring topology, reducing the communication probability $p$ does not lead to communication acceleration for ProxSkip \cite{ProxSkip}. This observation aligns with the experimental results reported in \cite{Alghunaim2023}. Nevertheless, the proposed $\textsc{MG-Skip}$ demonstrates a distinct advantage when it comes to local updates, as it enables communication acceleration under such conditions. From Fig. \ref{CASE1}, we also know that we can often skip the $\textsc{FastGoss}(\cdot)$ without any deterioration of the convergence rate.

\subsection{Decentralized Logistic Regression on ijcnn1 Dataset}

\begin{table*}[!t]
\renewcommand\arraystretch{1.6}
\caption{Experimental results of existing decentralized methods on ijcnn1 dataset. Stopping criterion: $\nicefrac{\|\M{x}^t - \M{x}^\star\|}{\|\M{x}^\star\|} < 10^{-7}$.}
\begin{center}
\scalebox{0.725}{
\begin{tabular}{|c|clcl|clcl|clcl|}
\hline
\multirow{2}{*}{Algorithm}  & \# of Iter. &  \multicolumn{1}{c}{speedup} & \# of Comm. & \multicolumn{1}{c|}{speedup}& \# of Iter. &  \multicolumn{1}{c}{speedup} & \# of Comm. & \multicolumn{1}{c|}{speedup}
& \# of Iter. &  \multicolumn{1}{c}{speedup} & \# of Comm. & \multicolumn{1}{c|}{speedup}   \\
\cline{2-13}
~& \multicolumn{4}{c|}{network settings: ring graph, $1-\rho=0.0576$}&\multicolumn{4}{c|}{network settings: $\iota=0.25$, $1-\rho=0.1186$}&\multicolumn{4}{c|}{network settings: $\iota=0.5$, $1-\rho=0.3857$}\\
\hline
Prox-GT \cite{Harnessing,Sulaiman2021} &$1820$&$\textcolor[rgb]{0.7,0,0}{\downarrow}$ $( 0.1407 \times)$&$1820\times 2$&$\textcolor[rgb]{0.7,0,0}{\downarrow}$ $(0.2110 \times)$&$518$&$\textcolor[rgb]{0.7,0,0}{\downarrow}$ $( 0.4942 \times)$&$518\times2$&$\textcolor[rgb]{0.7,0,0}{\downarrow}$ $( 0.4942 \times)$&$257$&$\textcolor[rgb]{0.7,0,0}{\downarrow}$ $( 0.9961 \times)$&$257\times2$&$\textcolor[rgb]{0.7,0,0}{\downarrow}$ $( 0.9961 \times)$\\
Prox-DIGing \cite{DIGing,Sulaiman2021}& $1966$&$\textcolor[rgb]{0.7,0,0}{\downarrow}$ $( 0.1302 \times)$&$1966$ $\times$ 2 &$\textcolor[rgb]{0.7,0,0}{\downarrow}$ $(0.3906 \times)$&$590$&$\textcolor[rgb]{0.7,0,0}{\downarrow}$ $( 0.4339 \times)$&590&$\textcolor[rgb]{0.7,0,0}{\downarrow}$ $( 0.4339 \times)$&$257$&$\textcolor[rgb]{0.7,0,0}{\downarrow}$ $( 0.9961 \times)$&$257\times2$&$\textcolor[rgb]{0.7,0,0}{\downarrow}$ $( 0.9961 \times)$\\
Prox-EXTRA \cite{EXTRA,Sulaiman2021} &$959$&$\textcolor[rgb]{0.7,0,0}{\downarrow}$ $( 0.2669 \times)$&$959$&$\textcolor[rgb]{0,0.6,0}{\uparrow}$ $(1.6017 \times)$&$483$&$\textcolor[rgb]{0.7,0,0}{\downarrow}$ $( 0.5300 \times)$&483&$\textcolor[rgb]{0,0.6,0}{\uparrow}$ $(1.0600 \times)$&$265$&$\textcolor[rgb]{0.7,0,0}{\downarrow}$ $( 0.9660 \times)$&$265$&$\textcolor[rgb]{0,0.6,0}{\uparrow}$ $(1.9320 \times)$\\
Prox-NIDS \cite{NIDS,Sulaiman2021} & $928$&$\textcolor[rgb]{0.7,0,0}{\downarrow}$ $( 0.2759 \times)$&$928$&$\textcolor[rgb]{0,0.6,0}{\uparrow}$ $(1.6552 \times)$&$463$&$\textcolor[rgb]{0.7,0,0}{\downarrow}$ $( 0.5529 \times)$
&$463$&$\textcolor[rgb]{0,0.6,0}{\uparrow}$ $(1.1058 \times)$&$256$&$\textcolor[rgb]{1.00,1.00,1.00}{\downarrow}$ $(1\times)$&$256$&$\textcolor[rgb]{0,0.6,0}{\uparrow}$ $(2 \times)$
\\
RandProx \cite{RandProx}&818 & $\textcolor[rgb]{0.7,0,0}{\downarrow}$ $(0.3130 \times)$& 655 & $\textcolor[rgb]{0,0.6,0}{\uparrow}$ $(2.3450 \times)$&$256$&$\textcolor[rgb]{1.00,1.00,1.00}{\downarrow}$ $(1\times)$&$205$&$\textcolor[rgb]{0,0.6,0}{\uparrow}$ $(2.4975 \times)$ &$256$&$\textcolor[rgb]{1.00,1.00,1.00}{\downarrow}$ $(1\times)$&$205$&$\textcolor[rgb]{0,0.6,0}{\uparrow}$ $(2.4975 \times)$ \\
\hline
\rowcolor[gray]{0.95}MG-SONATA \cite{SONATA}  &  &&&&  &&&&  &&&\\
\rowcolor[gray]{0.95}/ PMGT-VR \cite{PMGT-VR}  &  &&&&  &&&&  &&&\\
\rowcolor[gray]{0.95} \textbf{Baseline}
&\multirow{-3}{*}{$256$} & \multirow{-3}{*}{$\textcolor[gray]{0.95}{\downarrow}$ $(1 \times)$} &\multirow{-3}{*}{$768 \times 2$} &\multirow{-3}{*}{$\textcolor[gray]{0.95}{\downarrow}$ $(1 \times)$}
&\multirow{-3}{*}{$256$} & \multirow{-3}{*}{$\textcolor[gray]{0.95}{\downarrow}$ $(1 \times)$} &\multirow{-3}{*}{$256 \times 2$} &\multirow{-3}{*}{$\textcolor[gray]{0.95}{\downarrow}$ $(1 \times)$}
&\multirow{-3}{*}{$256$} & \multirow{-3}{*}{$\textcolor[gray]{0.95}{\downarrow}$ $(1 \times)$} &\multirow{-3}{*}{$256 \times 2$} &\multirow{-3}{*}{$\textcolor[gray]{0.95}{\downarrow}$ $(1 \times)$}
  \\
\hline
\rowcolor[gray]{0.9} $\textsc{MG-Skip}$, $p=1$&$256$& $\textcolor[gray]{0.9}{\downarrow}$ $(1\times)$ & $768$ & $\textcolor[rgb]{0,0.6,0}{\uparrow}$ $(2 \times)$&$256$&$\textcolor[gray]{0.9}{\downarrow}$ $(1\times)$&$256$&$\textcolor[rgb]{0,0.6,0}{\uparrow}$ $(2 \times)$&$256$&$\textcolor[gray]{0.9}{\downarrow}$ $(1\times)$&$256$&$\textcolor[rgb]{0,0.6,0}{\uparrow}$ $(2 \times)$ \\
\rowcolor[gray]{0.9} $\textsc{MG-Skip}$, $p=0.5$&$256$& $\textcolor[gray]{0.9}{\downarrow}$ $(1 \times)$  &$384$&$\textcolor[rgb]{0,0.6,0}{\uparrow}$ $(4 \times)$&$256$&$\textcolor[gray]{0.9}{\downarrow}$ $(1\times)$&$128$&$\textcolor[rgb]{0,0.6,0}{\uparrow}$ $(4 \times)$&$256$&$\textcolor[gray]{0.9}{\downarrow}$ $(1\times)$&$128$&$\textcolor[rgb]{0,0.6,0}{\uparrow}$ $(4 \times)$\\
\rowcolor[gray]{0.9} $\textsc{MG-Skip}$, $p=0.2$&$256$& $\textcolor[gray]{0.9}{\downarrow}$ $(1 \times)$ &$154$&$\textcolor[rgb]{0,0.6,0}{\uparrow}$ $(9.9740 \times)$&$256$ &$\textcolor[gray]{0.9}{\downarrow}$ $(1 \times)$&$52$
&$\textcolor[rgb]{0,0.6,0}{\uparrow}$ $(9.8461 \times)$&$256$ &$\textcolor[gray]{0.9}{\downarrow}$ $(1 \times)$&$52$
&$\textcolor[rgb]{0,0.6,0}{\uparrow}$ $(9.8461 \times)$
\\
\hline
\end{tabular}}
 \label{TableCOMPA-EX}
\end{center}
\end{table*}

\begin{figure*}[!t]
  \centering
  \setlength{\abovecaptionskip}{-2pt}
  \subfigure{
  \includegraphics[width=0.32\linewidth]{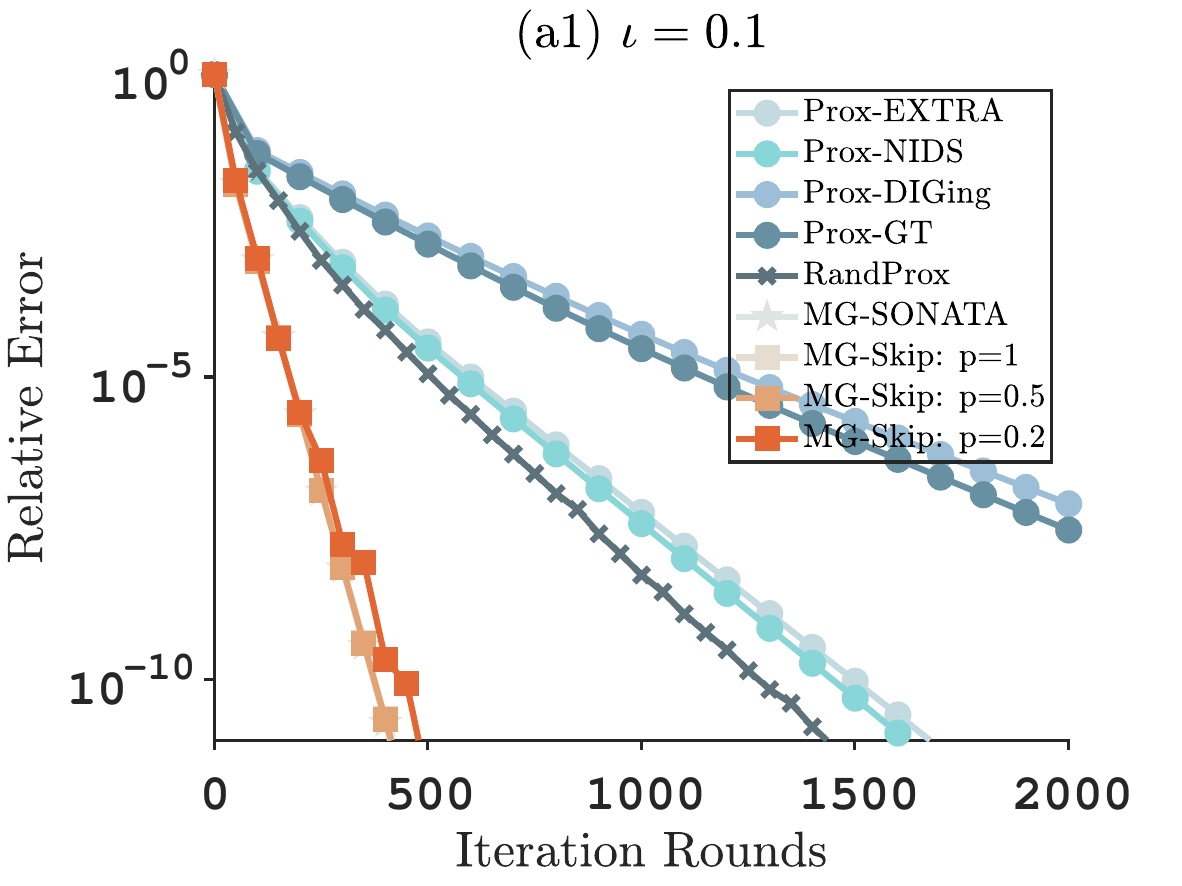}}
  \subfigure{
  \includegraphics[width=0.32\linewidth]{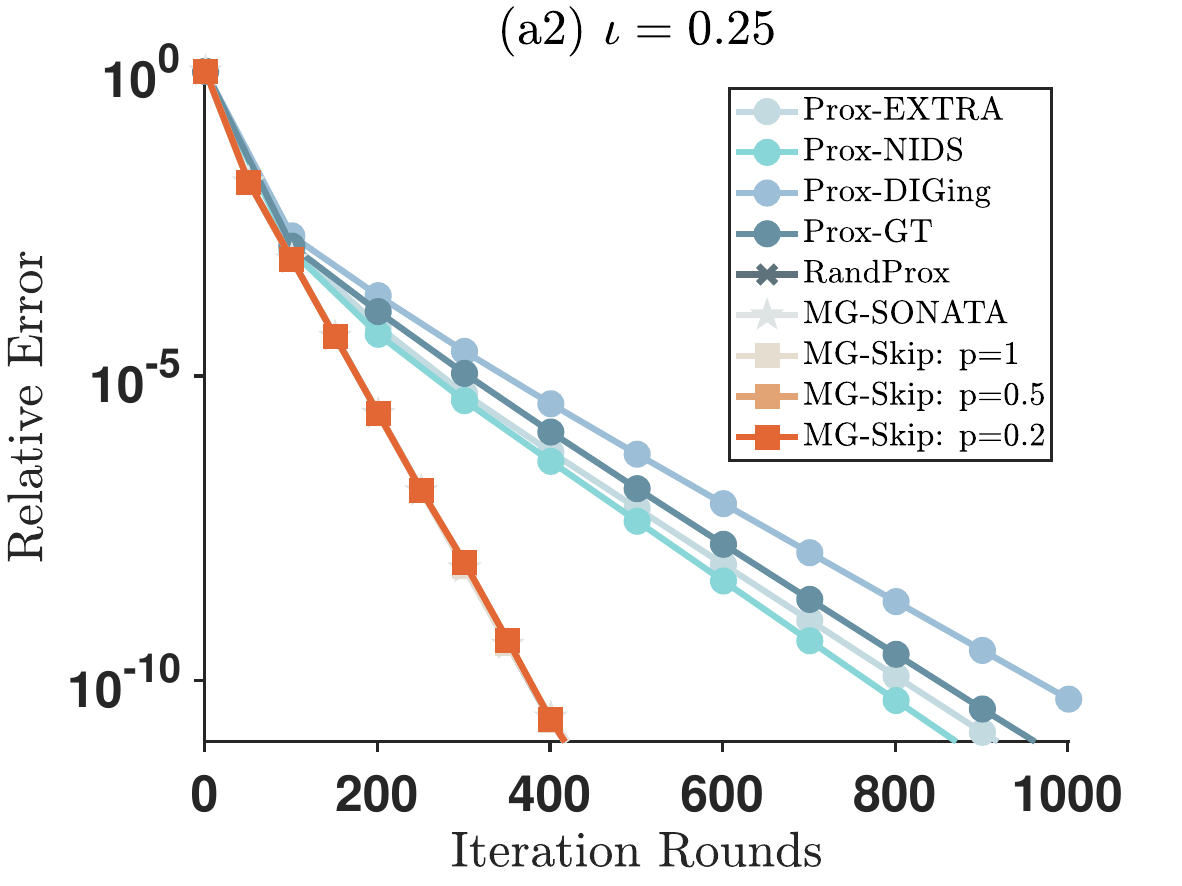}}
  \subfigure{
  \includegraphics[width=0.32\linewidth]{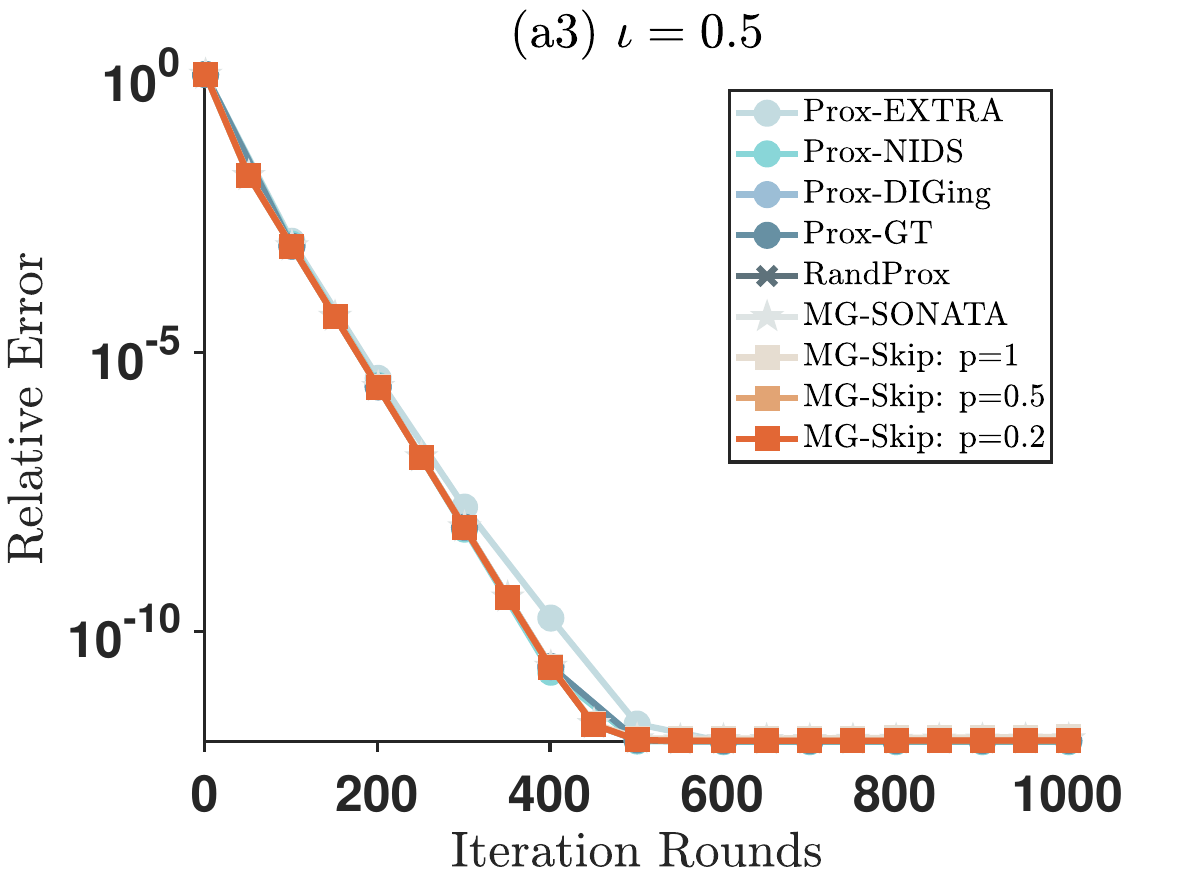}}
  \subfigure{
  \includegraphics[width=0.32\linewidth]{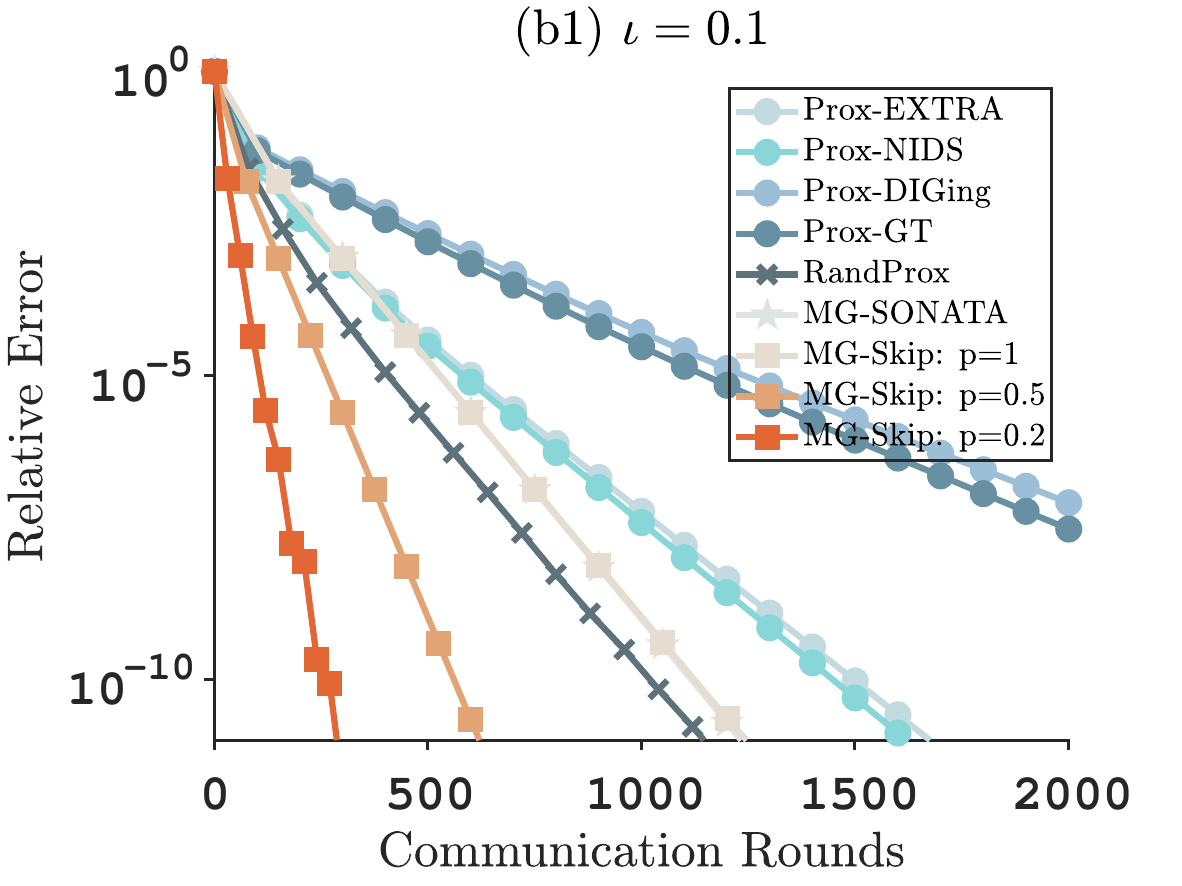}}
  \subfigure{
  \includegraphics[width=0.32\linewidth]{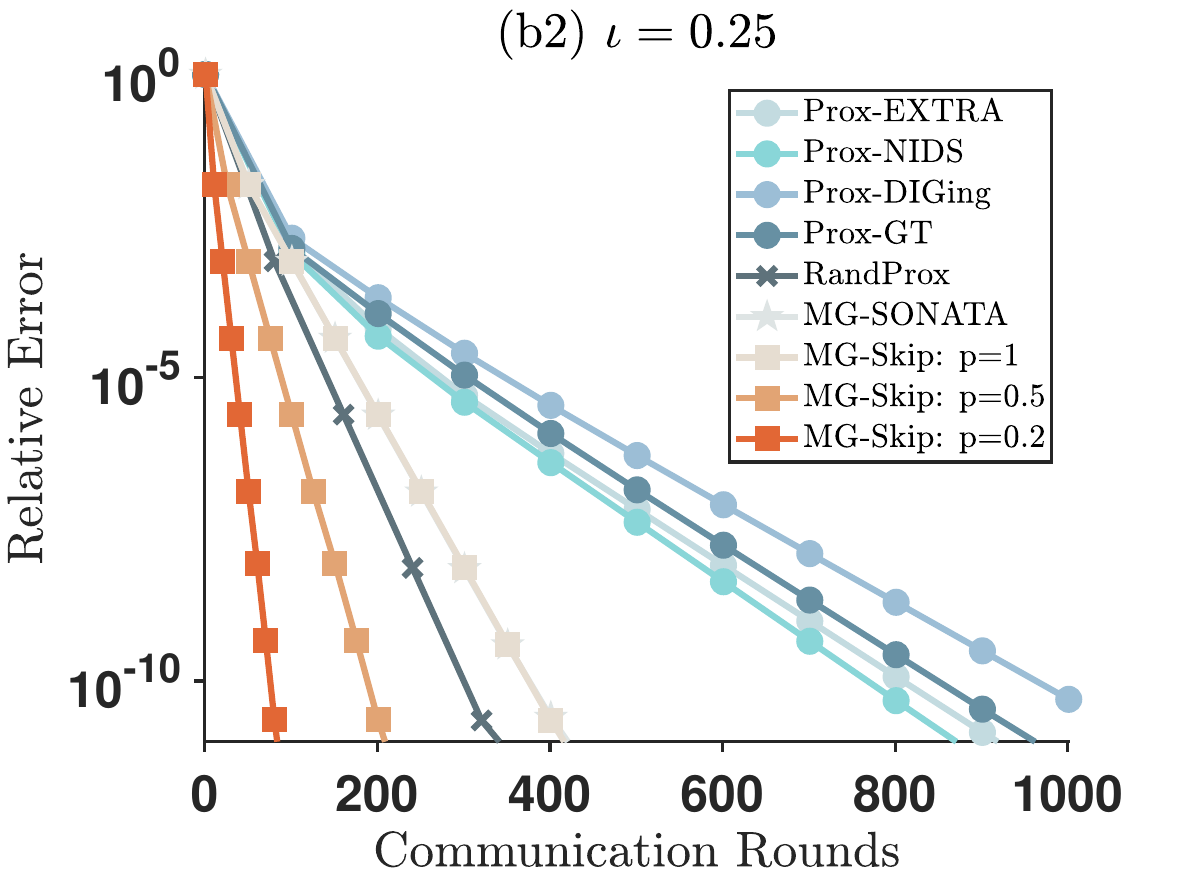}}
  \subfigure{
  \includegraphics[width=0.32\linewidth]{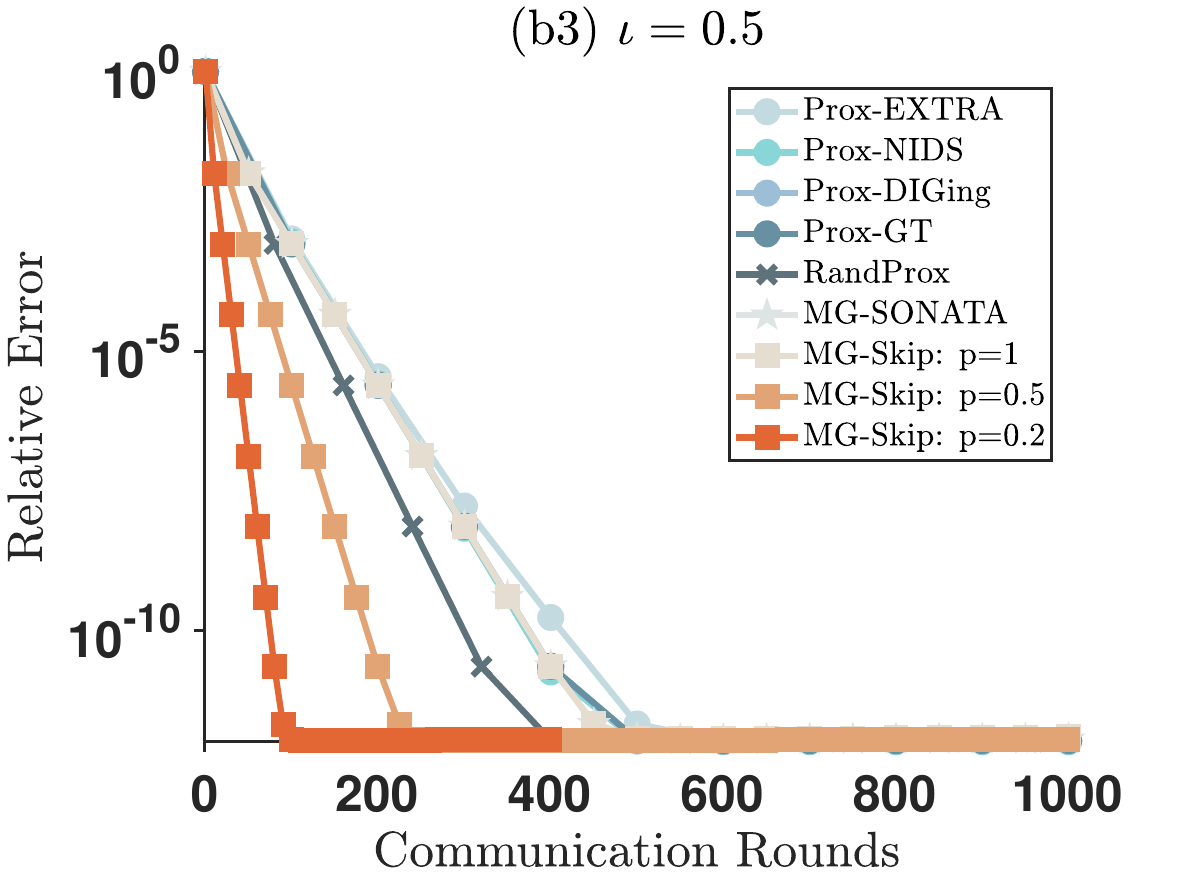}}
  \caption{Experimental comparison results for logistic regression problem with a nonsmooth regularizer $\|\m{x}\|_1$ on ijcnn1 dataset.}
  \label{CASE2}
\end{figure*}

In this subsection, a numerical experiment is given to demonstrate our findings on the decentralized logistic regression problem on ijcnn1 \cite{LibSVM} dataset. The loss function is defined as follows:
$$
h(\m{x})=\frac{1}{n}\sum_{i=1}^{n}\frac{1}{m_i}\sum_{j=1}^{m_{i}} \ln (1+e^{-(\mathcal{A}_{i j}\tr {\M{x}})\mathcal{B}_{i j}})+\gamma_1\|\m{x}\|^2+\gamma_2\|\m{x}\|_1.
$$
Here, any node $i$ holds its own training date $\left(\mathcal{A}_{i j}, \mathcal{B}_{i j}\right) \in$ $\mathbb{R}^{d} \times\{-1,1\}, j=1, \cdots, m_{i}$, including sample vectors $\mathcal{A}_{i j}$ and corresponding classes $\mathcal{B}_{i j}$, $\gamma_1=0.01$ and $\gamma_2=0.001$, $d=22$ and $\sum_{i=1}^{n}m_i=49950$. On this dataset, we have $L\thickapprox0.5$ and $\kappa\thickapprox25$. For all experiments, we first compute the solution $\m{x}^\star$ to \eqref{EQ:Problem1} by centralized methods, and then run over a randomly generated connected network with $N$ agents and $\frac{\iota n(n-1)}{2}$ undirected edges, where $\iota$ is the connectivity ratio. The mixing matrix $W$ is generated with the Metropolis-Hastings rule. In $\textsc{MG-Skip}$, we set the stepsize $\alpha=\nicefrac{1}{L}$. In other algorithms, we use the theoretical stepsizes provided in our respective literature.

The comparison results of communication and iteration rounds with Prox-EXTRA, Prox-NIDS, Prox-DIGing, Prox-GT, RandProx, and MG-SONATA are shown as Fig. \ref{CASE2}. To RandProx \cite{RandProx}, we set $p=\nicefrac{1}{\sqrt{(1-\rho)\kappa}}\approx0.8$. In Fig. \ref{CASE2} (b), we have that $\textsc{MG-Skip}$ with $p=\nicefrac{1}{\sqrt{\kappa}}\thickapprox0.2$ and $p=\nicefrac{2.5}{\sqrt{\kappa}}\thickapprox0.5$ outperform the other methods. In Fig. \ref{CASE2} (a), we have that $\textsc{MG-Skip}$ with $p=1,0.5,0.2$ have the same iteration performance. It implies that, when $p\geq\nicefrac{1}{\sqrt{\kappa}}$, skipping some communication will not affect the convergence rate. Although we note that RandProx achieves linear convergence in this numerical experiment, a provable linear convergence for decentralized composite optimization has not been reported. In addition, to further highlight the advantages of $\textsc{MG-Skip}$ compared to existing algorithms, we present the required iteration rounds and communication overhead in Table \ref{TableCOMPA-EX}, when $\nicefrac{\|\M{x}^t - \M{x}^\star\|}{\|\M{x}^\star\|} < 10^{-7}$. Here, ``$\times 2$'' indicates that two variables need to be communicated in each communication round. Setting $p = 0.2 \approx \nicefrac{1}{\sqrt{\kappa}}$, $\textsc{MG-Skip}$ achieves approximate $2\sqrt{\kappa}\approx10\times$ communication speedup compared to MG-SONATA. This improvement arises because the communication complexity of MG-SONATA and $\textsc{MG-Skip}$ are $\tilde{\mathcal{O}}(\nicefrac{\kappa}{\sqrt{1-\rho}})$ and $\tilde{\mathcal{O}}(\nicefrac{p\kappa}{\sqrt{1-\rho}})$, respectively, which are consistent with our theoretical findings and the comparative analysis presented in Table \ref{TableCOMPA}.

\begin{figure*}[!t]
  \centering
  \setlength{\abovecaptionskip}{-2pt}
  \subfigure{
  \includegraphics[width=0.32\linewidth]{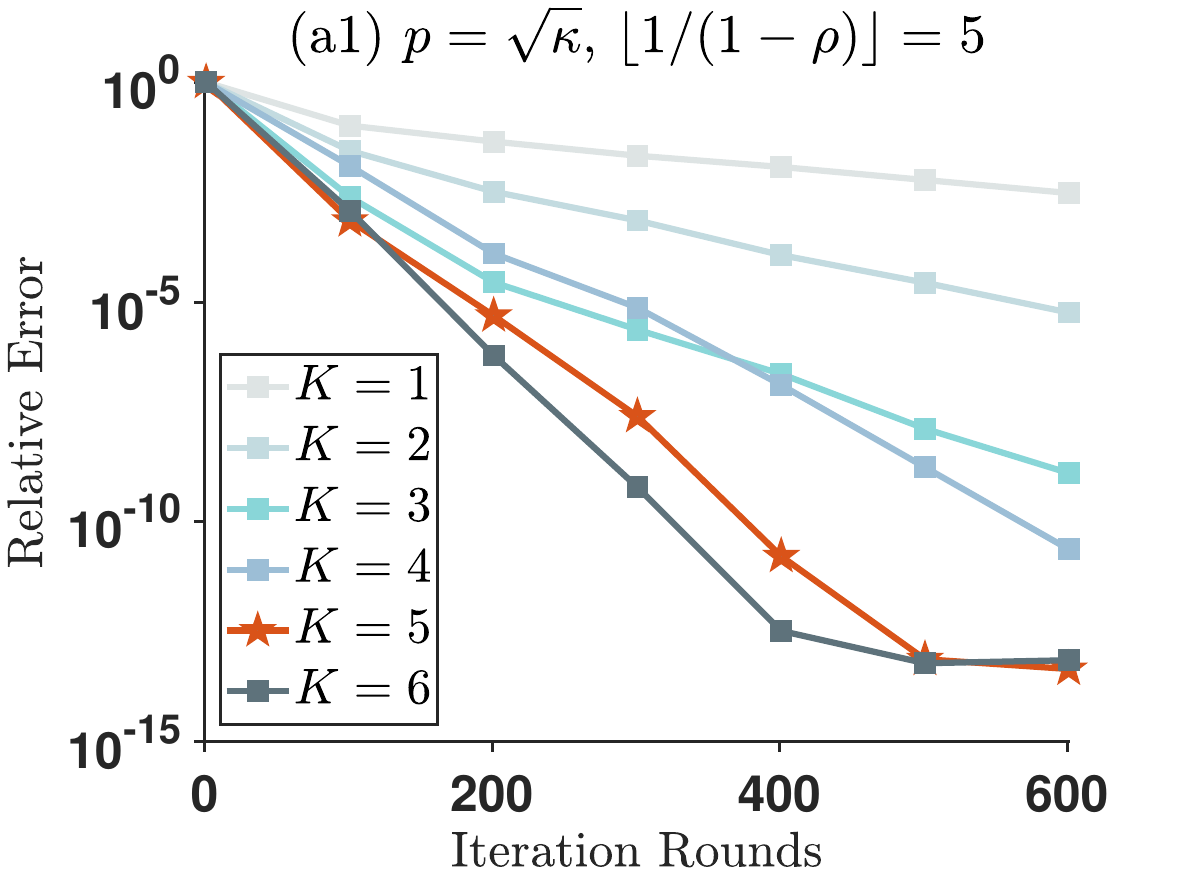}}
  \subfigure{
  \includegraphics[width=0.32\linewidth]{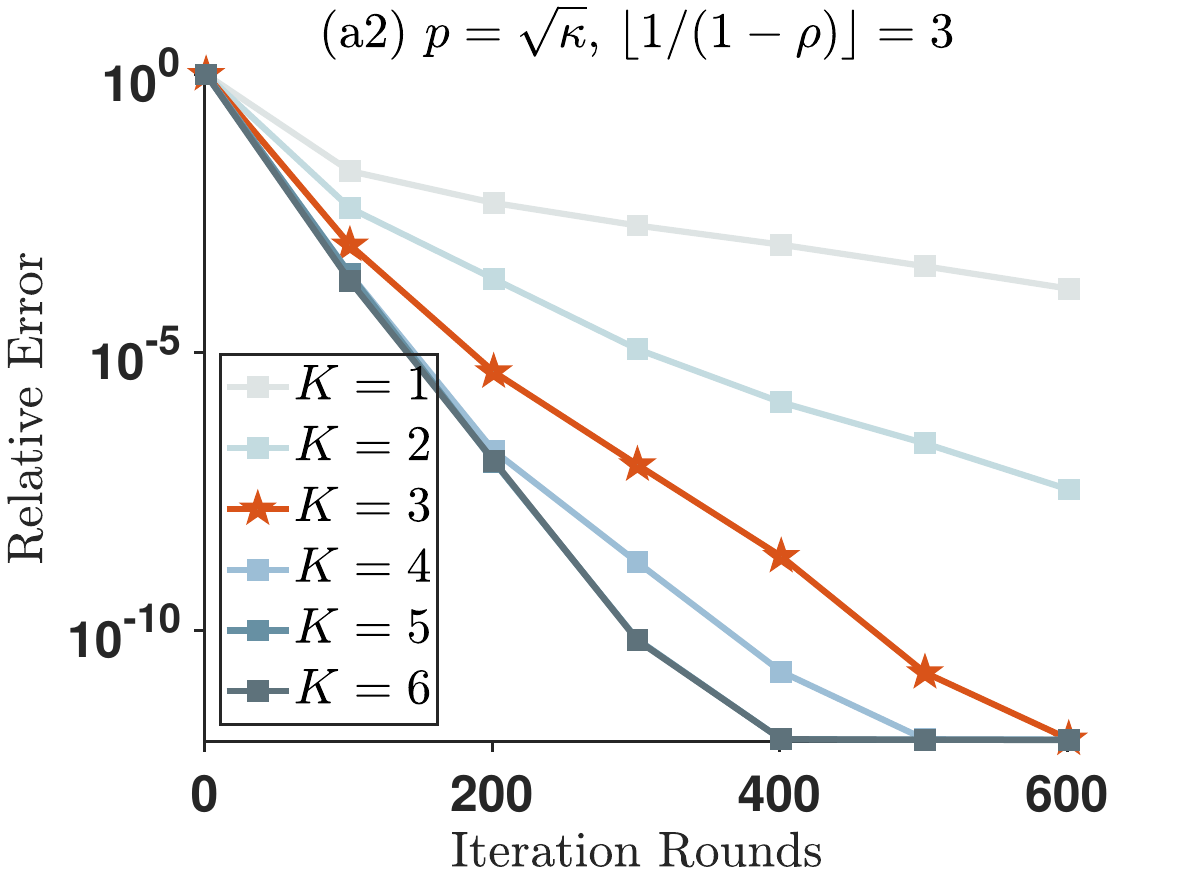}}
  \subfigure{
  \includegraphics[width=0.32\linewidth]{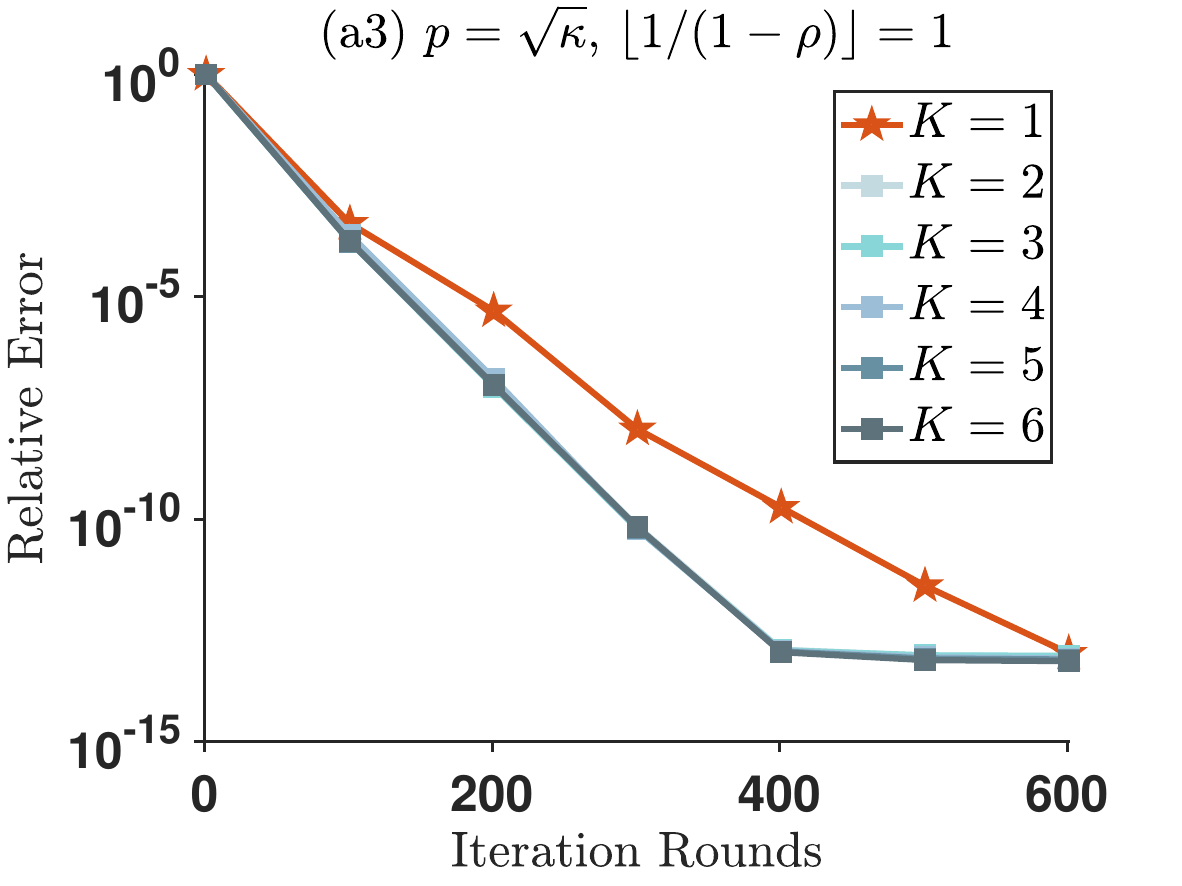}}
  \subfigure{
  \includegraphics[width=0.32\linewidth]{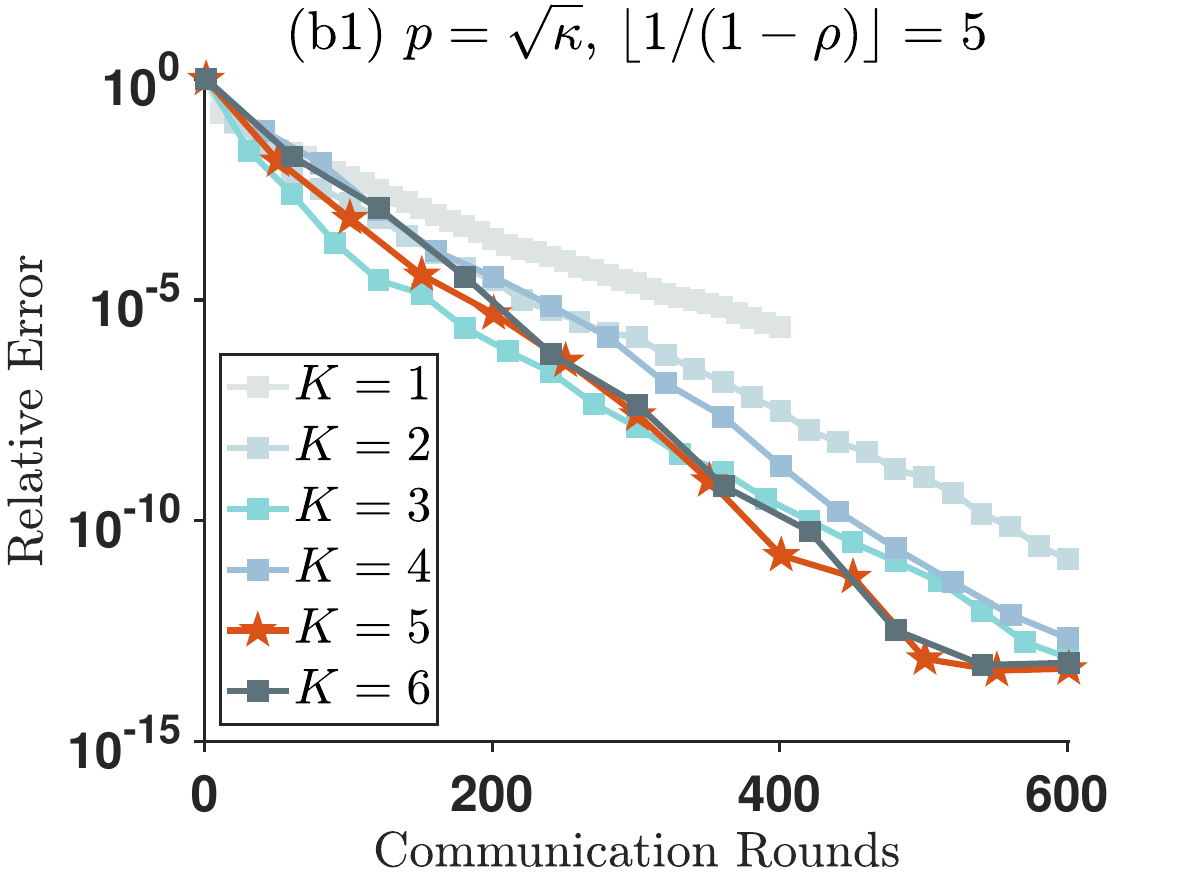}}
  \subfigure{
  \includegraphics[width=0.32\linewidth]{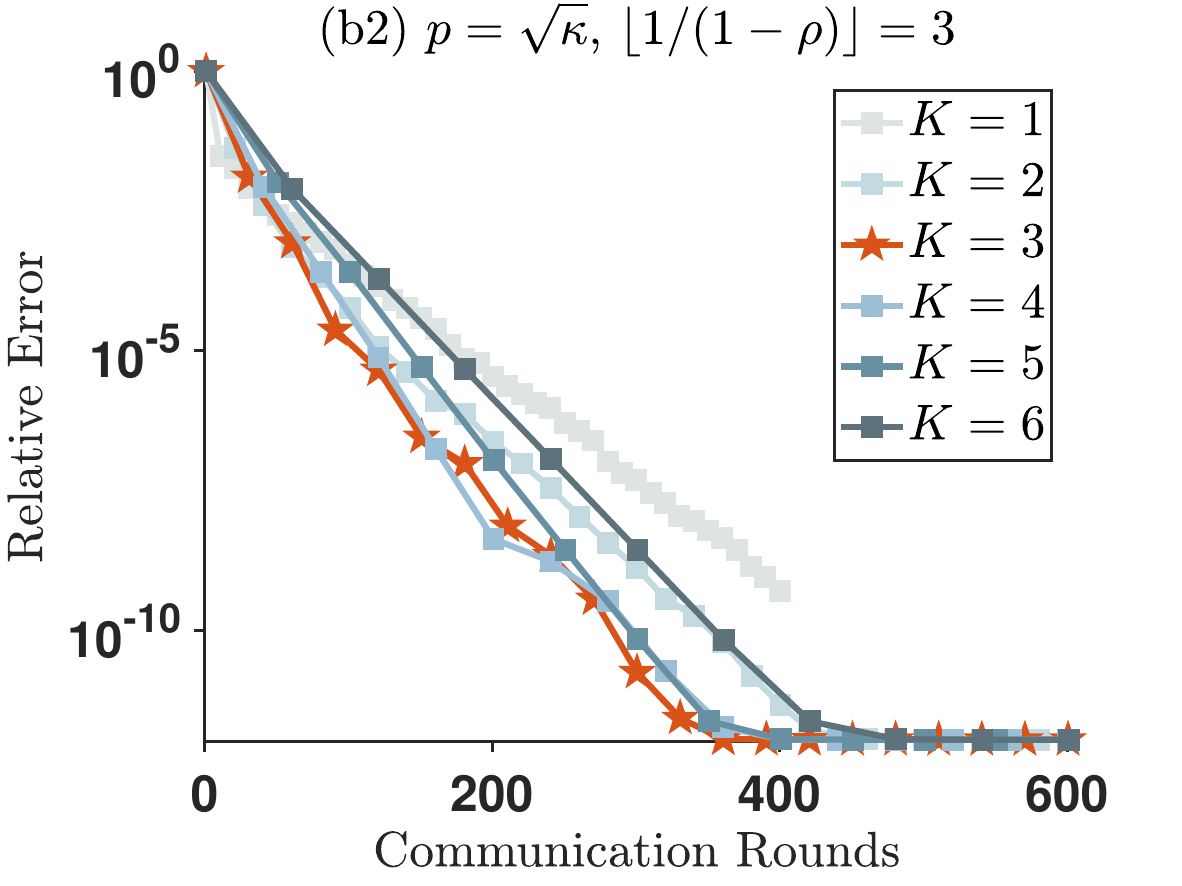}}
  \subfigure{
  \includegraphics[width=0.32\linewidth]{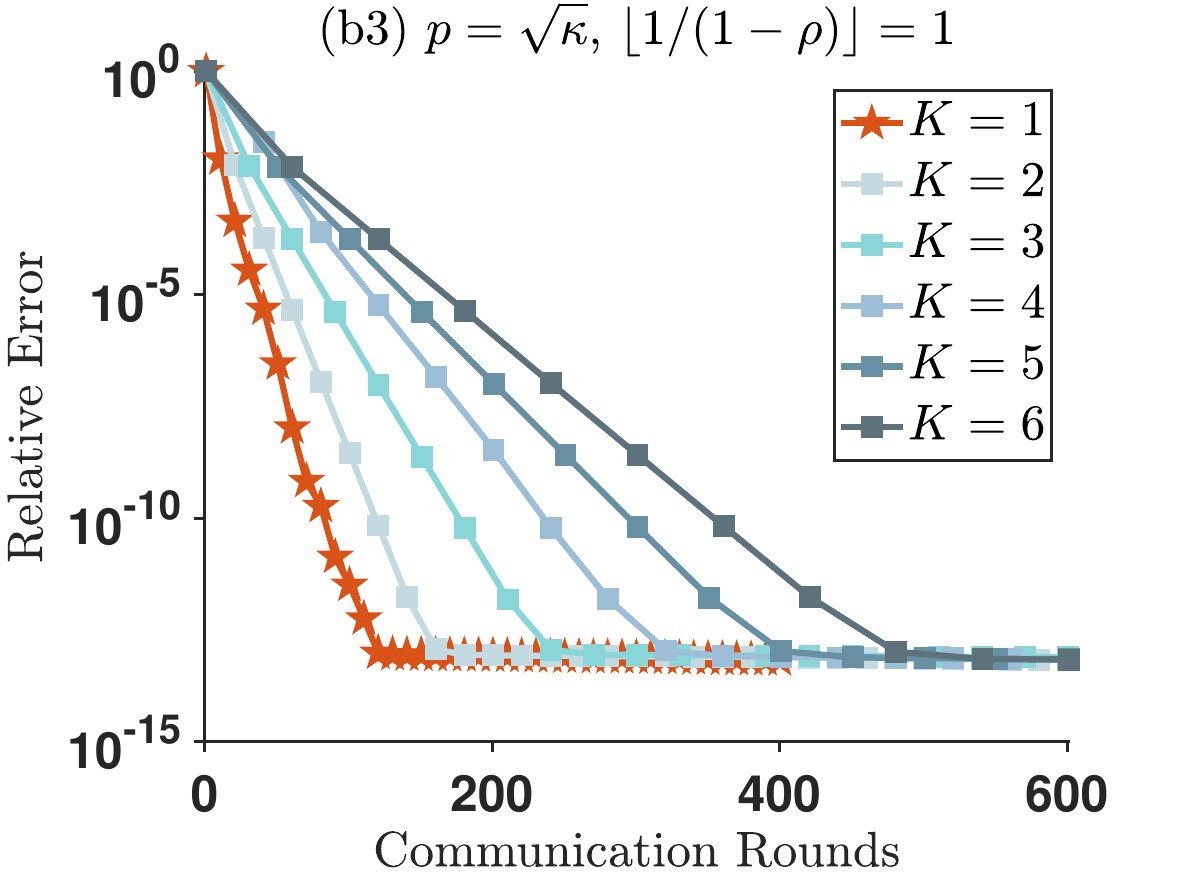}}
  \caption{Experimental results on the influence of consecutive communication rounds $K$ on the performance of $\textsc{MG-Skip}$ over different topologies.}
  \label{CASE2.2para}
\end{figure*}

Moreover, we continue our discussion on the influence of consecutive communication rounds $K$ on the performance of $\textsc{MG-Skip}$. According to the previous analysis, the theoretically optimal value of $p$ is $\nicefrac{1}{\sqrt{\kappa}}\approx0.2$. Thus, to all case, we set $p=0.2$. The experimental results are shown as Fig. \ref{CASE2.2para}. Experimental results show that multi-gossip communications can accelerate convergence. However, this acceleration is not endless, when a certain value of $K$ is reached, continuing to increase $K$ will not achieve acceleration. Through comparison, it can be concluded that when $K=\lfloor\nicefrac{1}{\sqrt{1-\rho}}\rfloor$, the communication efficiency of $\textsc{MG-Skip}$ is the highest (requiring the least total communication rounds), which is consistent with our theoretical results.

\begin{figure*}[!t]
  \centering
  \setlength{\abovecaptionskip}{-2pt}
  \subfigure{
  \includegraphics[width=0.3\linewidth]{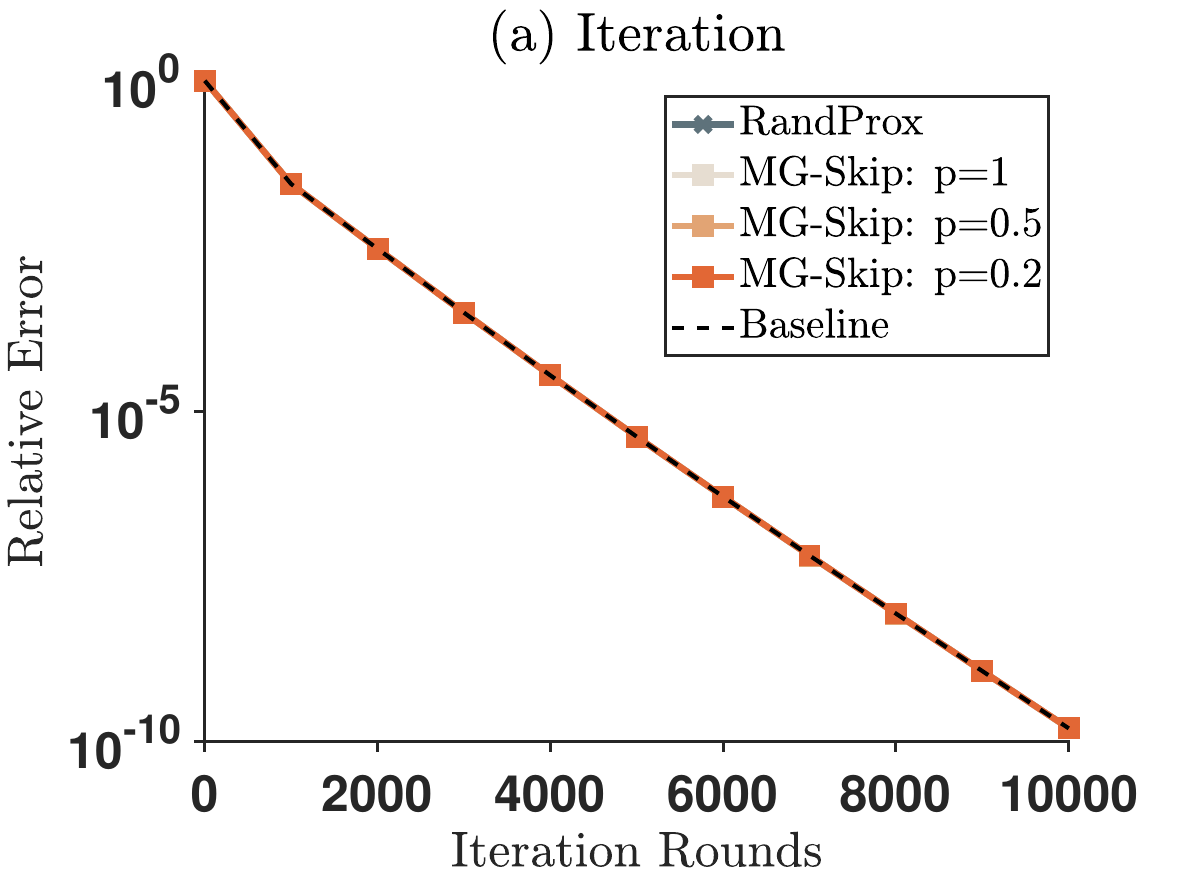}}
  \subfigure{
  \includegraphics[width=0.3\linewidth]{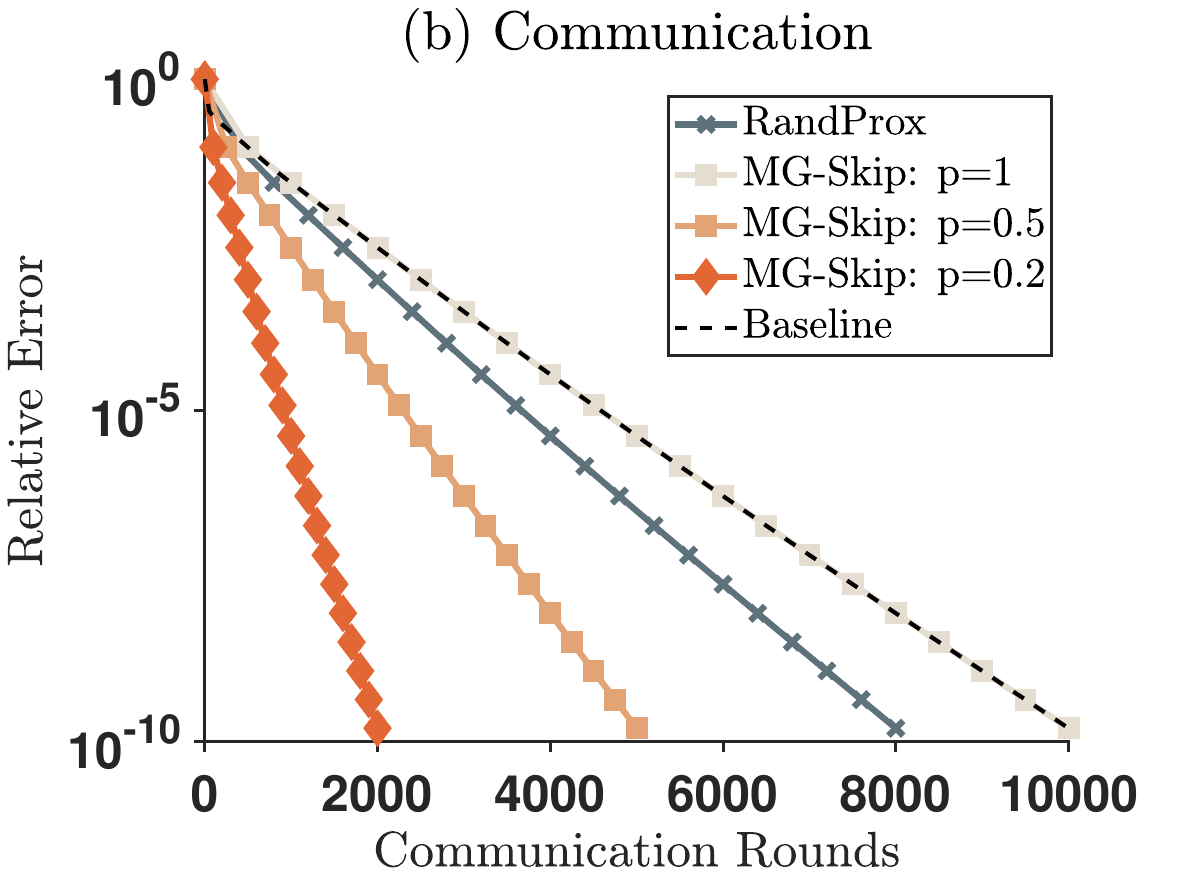}}
  \subfigure{
  \includegraphics[width=0.36\linewidth]{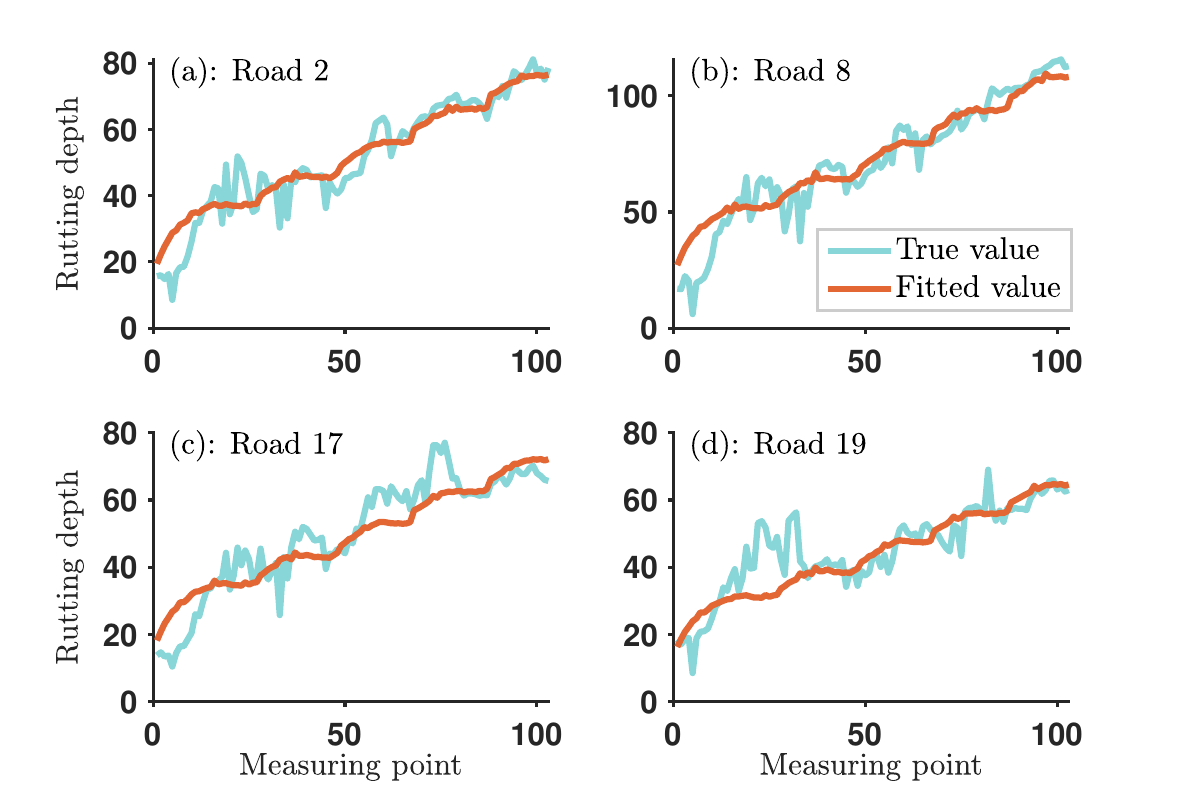}}
  \caption{Experimental results for sparse linear regression on rutting data. Here, baseline denotes MG-SONATA. The left image reports the relationship between the number of iterations and the relative error. The middle image reports the relationship between the number of communications and the relative error. The right image reports the differences between the fitted rutting depth ($\textsc{MG-Skip}:p=0.2$) and the actual rutting depth.}
  \label{CASE3.1}
\end{figure*}

\subsection{Decentralized Sparse Linear Regression on Rutting Data}
Rutting is an important index to study long-life asphalt pavement. In this subsection, we consider the decentralized sparse linear regression on rutting data. The experimental data used in this example was obtained from RIOHTrack (Research Institute of Highway MOT Track), and the structure of RIOHTrack is shown as \cite{RIOHTrack}. We choose the loss function as follows:
$$
h(\m{x})=\frac{1}{n}\sum_{i=1}^{n}\sum_{j=1}^{m_i}\frac{1}{2}\left\|\mathcal{A}_{ij}\m{x}-\mathcal{B}_{ij}\right\|^2
+\gamma_1\|\m{x}\|^2+\gamma_2\|\m{x}\|_1,
$$
where $\mathcal{A}_{ij}\in\mathbb{R}^{20}$ is the sample vector, $\mathcal{B}_{ij}\in\mathbb{R}_+$ is the rutting depth, $\gamma_1=0.01$ and $\gamma_2=0.001$. We used load, pavement temperature, and other data from Road 2, Road 8, Road 17, and Road 19 to fit their rutting in a decentralized manner. The communication topology is a ring topology with this 4 roads. Each road does not need to share private data, thus it provides better generalization and privacy.

The comparison results of communication rounds with RandProx and MG-SONATA are shown as Fig. \ref{CASE3.1}. In this numerical experiment, we take MG-SONATA as baseline. By comparison, we know that $\textsc{MG-Skip}$ outperform the other methods, when $p=0.5$ and $p=0.2$. Fig. \ref{CASE3.1} shows the differences between the fitted data ($\textsc{MG-Skip}:p=0.2$) and the actual data.

\section{Conclusion}\label{SEC6}
In this paper, we present $\textsc{MG-Skip}$, a novel algorithm that integrates probabilistic local updates for decentralized composite optimization. For the smooth (strongly convex)+nonsmooth (convex) composite form, we establish the linear convergence of $\textsc{MG-Skip}$, enabling it as the first algorithm with both random communication skipping and linear convergence, to the best of our knowledge.
Additionally our results not only demonstrate that $\textsc{MG-Skip}$ achieves the optimal communication complexity, but also confirm the benefits of local updates in the nonsmooth setup.
\appendices

\section{Proof of Proposition \ref{Pro1}}\label{AP:Pro1}
\begin{IEEEproof}
Since $\M{W}$ is symmetric and doubly stochastic, the matrix $\M{M}_K$ is symmetric and doubly stochastic.  From $\M{M}_{k+1}=(1+\eta)\M{W}\M{M}_k-\eta\M{M}_{k-1}$, we have
\begin{align*}
\left(
  \begin{array}{l}
     \M{M}_{k+1}\\
    \M{M}_{k} \\
  \end{array}
\right)=\underbrace{\left(
          \begin{array}{cc}
            (1+\eta)\M{W} & -\eta\M{I} \\
            \M{I} & 0 \\
          \end{array}
        \right)}_{\tilde{\M{W}}} \left(
  \begin{array}{l}
     \M{M}_{k}\\
    \M{M}_{k-1} \\
  \end{array}
\right).
\end{align*}
By \cite{Yuan2021}, the projection of the augmented matrix $\tilde{\M{W}}$ on the subspace orthogonal to $\M{1}$ is a contraction with spectral norm of $\frac{\rho}{1+\sqrt{1-\rho^2}}$.
Since
$$
\frac{\rho}{1+\sqrt{1-\rho^2}}\leq1-\sqrt{1-\rho}
$$
for any $\rho\in[0,1]$, we have
\begin{align*}
\left\|\left(
  \begin{array}{l}
     \M{M}_{k}\\
    \M{M}_{k-1} \\
  \end{array}
\right)\M{x}
\right\|_{\m{F}}&\leq \frac{\rho}{1+\sqrt{1-\rho^2}}
\left\|\left(
  \begin{array}{l}
     \M{M}_{k-1}\\
    \M{M}_{k-2} \\
  \end{array}
\right)\M{x}
\right\|_{\m{F}}\\
&\leq\Big(1-\sqrt{1-\rho}\Big)^k\left\|\left(
  \begin{array}{l}
     \M{M}_{0}\\
    \M{M}_{-1} \\
  \end{array}
\right)\M{x}
\right\|_{\m{F}}\\
&=\sqrt{2}\Big(1-\sqrt{1-\rho}\Big)^k\Big\|\M{x}\Big\|_{\m{F}}
\end{align*}
for any $\M{x}\perp \M{1}_n$, i.e., $\M{1}_n\tr\M{x}=0$, where $\M{x}\in \mathbb{R}^{n\times d}$. Noting that $\bar{\M{M}}=\M{M}_K$, for any $\M{x}\perp \M{1}_n$, we have
$$
\Big\|\bar{\M{M}}\M{x}\Big\|_{\m{F}}\leq\sqrt{2}\Big(1-\sqrt{1-\rho}\Big)^K\Big\|\M{x}\Big\|_{\m{F}}.
$$
Since $\M{M}_K$ is symmetric and doubly stochastic, it holds that
$$
\varpi(\bar{\M{M}}-\nicefrac{1}{n}\M{11}\tr)\leq\sqrt{2}\Big(1-\sqrt{1-\rho}\Big)^K.
$$
Setting $K=\lfloor\frac{1}{\sqrt{1-\rho}}\rfloor$, we have that
$$
\varpi(\bar{\M{M}}-\nicefrac{1}{n}\M{11}\tr)\leq\sqrt{2}\Big(1-\sqrt{1-\rho}\Big)^{\frac{1}{\sqrt{1-\rho}}}\leq\frac{\sqrt{2}}{e},
$$
where $e$ is the Euler's number. Note that $\bar{\M{M}}$ is symmetric and doubly stochastic. The largest eigenvalue of the mixing matrix $\bar{\M{M}}$ is $\bar{\lambda}_1=1$, and the remaining eigenvalues are
$$
1>\bar{\lambda}_2\geq\bar{\lambda}_3\geq\cdots\geq\bar{\lambda}_n>-1.
$$
and
$$
\varpi(\bar{\M{M}}-\nicefrac{1}{n}\M{11}\tr)=\max\{|\bar{\lambda}_2|,|\bar{\lambda}_n|\}\in(0,1).
$$
Therefore, it holds that
$$
\sigma_m(\M{I}-\bar{\M{M}})\geq1-\varpi(\M{M}^K-\nicefrac{1}{n}\M{11}\tr)\geq1-\frac{\sqrt{2}}{e}\geq\frac{2}{5},
$$
where $\sigma_m(\M{I}-\bar{\M{M}})$ is the smallest none zero eigenvalue of $\M{I}-\bar{\M{M}}$.
\end{IEEEproof}

\section{Proof of Proposition \ref{Pro2}}\label{AP:Pro2}
\begin{IEEEproof}
Firstly, we assume $(\M{x}^\star,\M{z}^\star,\M{u}^\star)$ satisfying that $0\in\mathcal{M}(\M{x}^\star,\M{z}^\star,\M{u}^\star)$, i.e.,
\begin{align*}
0\in
\left(
                                     \begin{array}{c}
\M{z}^\star-\M{x}^\star+\alpha \nabla F(\M{x}^\star)+\alpha \sqrt{\nicefrac{1}{2}(\M{I}-\bar{\M{M}})} \M{u}^\star\\
-\sqrt{\nicefrac{1}{2}(\M{I}-\bar{\M{M}})}\M{z}^\star\\
(\M{x}^\star+\alpha \partial R(\M{x}^\star))-\left(\M{I}-{\nicefrac{1}{2}(\M{I}-\bar{\M{M}})}\right)\M{z}^\star
                                     \end{array}
                                   \right).
\end{align*}
Since $R(\M{x})=\sum_{i=1}^n r(\m{x}_i)$, it gives that
$$
\m{z}_i=\m{x}_i+\alpha \partial r(\m{x}_i)\Longrightarrow \m{x}_i=(I+\alpha\partial r)^{-1}(\m{z}_i),
$$
where $I$ is the identity mapping. Therefore, when $\sqrt{\nicefrac{1}{2}(\M{I}-\bar{\M{M}})}\M{z}^\star=0$, it gives that $\sqrt{\nicefrac{1}{2}(\M{I}-\bar{\M{M}})}\M{x}^\star=0$. Then, by $\M{z}^\star-\M{x}^\star=\alpha \partial R(\M{x}^\star)$, we have
\begin{align*}
0\in\left(
         \begin{array}{cc}
           \nabla F(\M{x}^\star)+\partial R(\M{x}^\star)+\sqrt{\nicefrac{1}{2}(\M{I}-\bar{\M{M}})}\M{u}^\star \\
           -\sqrt{\nicefrac{1}{2}(\M{I}-\bar{\M{M}})}\M{x}^\star \\
         \end{array}
       \right),
\end{align*}
i.e., $0\in\partial L(\M{x}^\star,\M{u}^\star)$.

Then, we assume $(\M{x}^\star,\M{u}^\star)$ satisfying that $0\in\partial L(\M{x}^\star,\M{u}^\star)$ and $\M{z}^\star=\M{x}^\star+\alpha \partial R(\M{x}^\star)$. In this case, it hold that
\begin{align*}
0\in\left(
         \begin{array}{cc}
           -\sqrt{\nicefrac{1}{2}(\M{I}-\bar{\M{M}})}\M{x}^\star \\
           \M{x}^\star+\alpha \partial R(\M{x}^\star)-\M{z}^\star
         \end{array}
       \right).
\end{align*}
Noting that $\m{z}_i=\m{x}_i+\alpha \partial r(\m{x}_i)$, we have $\sqrt{\nicefrac{1}{2}(\M{I}-\bar{\M{M}})}\M{x}^\star=0$ implies that $\sqrt{\nicefrac{1}{2}(\M{I}-\bar{\M{M}})}\M{z}^\star=0$. Then, we have
\begin{align*}
0\in\left(
         \begin{array}{cc}
           -\sqrt{\nicefrac{1}{2}(\M{I}-\bar{\M{M}})}\M{z}^\star \\
           \M{x}^\star+\alpha \partial R(\M{x}^\star)-\left(\M{I}-{\nicefrac{1}{2}(\M{I}-\bar{\M{M}})}\right)\M{z}^\star
         \end{array}
       \right).
\end{align*}
Therefore, $0\in\mathcal{M}(\M{x}^\star,\M{z}^\star,\M{u}^\star)$.
\end{IEEEproof}

\section{Proof of Lemma \ref{LEM-1}}\label{AP:LEM-1}
\begin{IEEEproof}
It follows from \eqref{KKT-Condition-2}, i.e., $\sqrt{\nicefrac{1}{2}(\M{I}-\bar{\M{M}})}\M{z}^\star=0$, and Proposition \ref{Pro1} that
the block elements of $\M{z}^\star$ are equal to each other, i.e.,
$\M{z}^\star=[\m{z}^\star,\ldots,\m{z}^\star]\tr$.
Let $\M{x}^\star=[\m{x}_1^\star,\ldots,\m{x}_n^\star]\tr$. By \eqref{KKT-Condition-3} and the definition of the proximal operator, we know
$$
\M{x}^{\star}=\mathrm{prox}_{\alpha R}(\M{z}^{\star})=\mathrm{argmin}_{\M{x}\in\mathbb{R}^{n\times d}}\{\alpha R(\M{x})+\nicefrac{1}{2}\|\M{x}-\M{z}^\star\|_{\mathrm{F}}^2\},
$$
which implies that
\begin{align}\label{proof-lem-eq1}
\m{x}_i^{\star}=\mathrm{argmin}_{\m{x}\in\mathbb{R}^{d}}\{\alpha r(\m{x})+\nicefrac{1}{2}\|\m{x}-\m{z}^\star\|^2\}.
\end{align}
Thus, we have $\m{x}_1^\star=\ldots=\m{x}_n^\star\triangleq \m{x}^\star$. In addition, from \eqref{proof-lem-eq1}, it holds that
\begin{align}\label{proof-lem-eq2}
0\in\partial r(\m{x}^\star)+\nicefrac{1}{\alpha}(\m{x}^\star-\m{z}^\star).
\end{align}
Note that $\M{x}^\star=[\m{x}^\star,\ldots,\m{x}^\star]\tr$, $\M{z}^\star=[\m{z}^\star,\ldots,\m{z}^\star]\tr$, and $\M{1}\tr\sqrt{\nicefrac{1}{2}(\M{I}-\bar{\M{M}})}=0$. Multiplying $\M{1}\tr$ from the left to both sides of equation \eqref{KKT-Condition-1}, it gives that
\begin{align}\label{proof-lem-eq3}
&n\m{z}^{\star}=n\m{x}^{\star}-\alpha \sum_{i=1}^{n}\nabla f_i(\m{x}^{\star})\nonumber\\
&\Rightarrow \nicefrac{1}{\alpha}(\m{x}^{\star}-\m{z}^{\star})=\frac{1}{n}\sum_{i=1}^{n}\nabla f_i(\m{x}^{\star}).
\end{align}
Combining \eqref{proof-lem-eq2} and \eqref{proof-lem-eq3}, we have $0\in\partial r(\m{x}^\star)+\frac{1}{n}\sum_{i=1}^{n}\nabla f_i(\m{x}^{\star})$.
Thus, $\m{x}^\star$ is a solution to problem \eqref{EQ:Problem1}.
\end{IEEEproof}

\section{Proof of Lemma \ref{LEM-2}}\label{AP:LEM-2}
\begin{IEEEproof}
We first define the following notations
\begin{align*}
\M{v}^t&:=\M{x}^t-\alpha \nabla F(\M{x}^t),\\
\M{v}^{\star}&:=\M{x}^{\star}-\alpha \nabla F(\M{x}^{\star}),\\
\hat{\M{u}}^{t+1}&:=\M{u}^t+\frac{p}{\alpha} \sqrt{\nicefrac{1}{2}(\M{I}-\bar{\M{M}})}{\M{z}}^{t},\\
\hat{\M{z}}^{t+1}&:={\M{z}}^t-\frac{\theta_t\alpha}{p}\sqrt{\nicefrac{1}{2}(\M{I}-\bar{\M{M}})}\left(\hat{\M{u}}^{t+1}-\M{u}^t\right).
\end{align*}

Taking the expectation over $\theta_t$, we have
\begin{align}\label{PROOF-IEQ1}
&\mathbb{E} \left[ \left\|\hat{\M{z}}^{t+1}-\M{z}^{\star} \right\|_{\m{F}}^2 ~|~ \theta_t \right]\nonumber\\
&=p\Big\|({\M{z}}^t-\M{z}^{\star})-\frac{\alpha}{p}\sqrt{\nicefrac{1}{2}(\M{I}-\bar{\M{M}})}\left(\hat{\M{u}}^{t+1}-\M{u}^t\right)\Big\|_{\m{F}}^2\nonumber\\
&\quad+(1-p)\Big\|{\M{z}}^t-\M{z}^{\star}\Big\|_{\m{F}}^2\nonumber\\
&= \Big\|{\M{z}}^t-\M{z}^{\star} \Big\|_{\m{F}}^2+\frac{\alpha^2}{p} \Big\|\sqrt{\nicefrac{1}{2}(\M{I}-\bar{\M{M}})}(\hat{\M{u}}^{t+1}-\M{u}^t) \Big\|_{\m{F}}^2\nonumber\\
&\quad-2\alpha \Big\langle{\M{z}}^t-\M{z}^{\star},\sqrt{\nicefrac{1}{2}(\M{I}-\bar{\M{M}})}(\hat{\M{u}}^{t+1}-\M{u}^t) \Big\rangle.
\end{align}
Since
\begin{align*}
\M{z}^t&=\M{x}^t-\alpha \nabla F(\M{x}^t)-\alpha\sqrt{\nicefrac{1}{2}(\M{I}-\bar{\M{M}})}\M{u}^t\\
&=\M{v}^t-\alpha\sqrt{\nicefrac{1}{2}(\M{I}-\bar{\M{M}})}\M{u}^t,
\end{align*}
and
\begin{align*}
\M{z}^\star&=\M{x}^\star-\alpha \nabla F(\M{x}^\star)-\alpha\sqrt{\nicefrac{1}{2}(\M{I}-\bar{\M{M}})}\M{u}^\star\\
&=\M{v}^\star-\alpha\sqrt{\nicefrac{1}{2}(\M{I}-\bar{\M{M}})}\M{u}^\star,
\end{align*}
it gives that
$$
{\M{z}}^t-\M{z}^{\star} = \M{v}^t-\M{v}^\star-\alpha\sqrt{\nicefrac{1}{2}(\M{I}-\bar{\M{M}})}(\M{u}^t-\M{u}^\star).
$$
Thus, we have
\begin{align}\label{PROOF-IEQ2}
&\left\|{\M{z}}^t-\M{z}^{\star} \right\|_{\m{F}}^2\nonumber\\
&=\Big\langle\M{v}^t-\M{v}^\star-\alpha\sqrt{\nicefrac{1}{2}(\M{I}-\bar{\M{M}})}(\M{u}^t-\M{u}^\star),{\M{z}}^t-\M{z}^{\star}\Big\rangle.
\end{align}
Combining \eqref{PROOF-IEQ1} with \eqref{PROOF-IEQ2}, it holds that
\begin{align}\label{PROOF-IEQ3}
&\mathbb{E} \left[ \left\|\hat{\M{z}}^{t+1}-\M{z}^{\star} \right\|_{\m{F}}^2 ~|~\theta_t\right]\nonumber\\
&=\underbrace{\Big\langle\M{v}^t-\M{v}^\star-\alpha\sqrt{\nicefrac{1}{2}(\M{I}-\bar{\M{M}})}(\M{u}^t-\M{u}^\star),{\M{z}}^t-\M{z}^{\star}\Big\rangle}_{\mathrm{i}}\nonumber\\
&\quad +\frac{\alpha^2}{p} \Big\|\sqrt{\nicefrac{1}{2}(\M{I}-\bar{\M{M}})}(\hat{\M{u}}^{t+1}-\M{u}^t) \Big\|_{\m{F}}^2\nonumber\\
&\quad\underbrace{-2\alpha \Big\langle{\M{z}}^t-\M{z}^{\star},\sqrt{\nicefrac{1}{2}(\M{I}-\bar{\M{M}})}(\hat{\M{u}}^{t+1}-\M{u}^t) \Big\rangle}_{\mathrm{ii}}.
\end{align}
Noting that
\begin{align*}
\mathrm{i}&=2\alpha \Big\langle\sqrt{\nicefrac{1}{2}(\M{I}-\bar{\M{M}})}(\M{u}^t-\M{u}^{\star}),\M{z}^t-\M{z}^{\star}\Big\rangle\nonumber\\
&\quad-2\alpha \Big\langle\sqrt{\nicefrac{1}{2}(\M{I}-\bar{\M{M}})}(\hat{\M{u}}^{t+1}-\M{u}^{\star}),\M{z}^t-\M{z}^{\star}\Big\rangle,
\end{align*}
it gives that
\begin{align}\label{PROOF-IEQ3-NEW1}
\mathrm{i}+\mathrm{ii}&=\Big\langle\M{v}^t-\M{v}^\star-\alpha\sqrt{\nicefrac{1}{2}(\M{I}-\bar{\M{M}})}(\M{u}^t-\M{u}^\star),\M{z}^t-\M{z}^{\star}\Big\rangle\nonumber\\
&\quad+2\alpha \Big\langle\sqrt{\nicefrac{1}{2}(\M{I}-\bar{\M{M}})}(\M{u}^t-\M{u}^{\star}),\M{z}^t-\M{z}^{\star}\Big\rangle\nonumber\\
&\quad-2\alpha \Big\langle\sqrt{\nicefrac{1}{2}(\M{I}-\bar{\M{M}})}(\hat{\M{u}}^{t+1}-\M{u}^{\star}),\M{z}^t-\M{z}^{\star}\Big\rangle\nonumber\\
&=\Big\langle\M{v}^t-\M{v}^\star+\alpha\sqrt{\nicefrac{1}{2}(\M{I}-\bar{\M{M}})}(\M{u}^t-\M{u}^\star),\M{z}^t-\M{z}^{\star}\Big\rangle\nonumber\\
&\quad-2\alpha \Big\langle\sqrt{\nicefrac{1}{2}(\M{I}-\bar{\M{M}})}(\hat{\M{u}}^{t+1}-\M{u}^{\star}),\M{z}^t-\M{z}^{\star}\Big\rangle.
\end{align}
Combining \eqref{PROOF-IEQ3} and \eqref{PROOF-IEQ3-NEW1}, it deduces that
\begin{align}\label{PROOF-IEQ3-NEW2}
&\mathbb{E} \left[ \left\|\hat{\M{z}}^{t+1}-\M{z}^{\star} \right\|_{\m{F}}^2 ~|~\theta_t\right]\nonumber\\
&=\Big\langle\M{v}^t-\M{v}^\star+\alpha\sqrt{\nicefrac{1}{2}(\M{I}-\bar{\M{M}})}(\M{u}^t-\M{u}^\star),\M{z}^t-\M{z}^{\star}\Big\rangle\nonumber\\
&\quad+\frac{\alpha^2}{p} \Big\|\sqrt{\nicefrac{1}{2}(\M{I}-\bar{\M{M}})}(\hat{\M{u}}^{t+1}-\M{u}^t) \Big\|_{\m{F}}^2\nonumber\\
&\quad-2\alpha \Big\langle\sqrt{\nicefrac{1}{2}(\M{I}-\bar{\M{M}})}(\hat{\M{u}}^{t+1}-\M{u}^{\star}),\M{z}^t-\M{z}^{\star}\Big\rangle.
\end{align}
Since $$
{\M{z}}^t-\M{z}^{\star} = \M{v}^t-\M{v}^\star-\alpha\sqrt{\nicefrac{1}{2}(\M{I}-\bar{\M{M}})}(\M{u}^t-\M{u}^\star),
$$
it holds that
\begin{align}\label{PROOF-IEQ4}
&\Big\langle\M{v}^t-\M{v}^\star+\alpha\sqrt{\nicefrac{1}{2}(\M{I}-\bar{\M{M}})}(\M{u}^t-\M{u}^\star),\M{z}^t-\M{z}^{\star}\Big\rangle\nonumber\\
&=\Big\langle(\M{v}^t-\M{v}^\star)+\alpha\sqrt{\nicefrac{1}{2}(\M{I}-\bar{\M{M}})}(\M{u}^t-\M{u}^\star),\nonumber\\
&\quad\quad(\M{v}^t-\M{v}^\star)-\alpha\sqrt{\nicefrac{1}{2}(\M{I}-\bar{\M{M}})}(\M{u}^t-\M{u}^\star)\Big\rangle\nonumber\\
&=\Big\|\M{v}^t-\M{v}^\star\Big\|_{\m{F}}^2-\alpha^2\Big\|\sqrt{\nicefrac{1}{2}(\M{I}-\bar{\M{M}})}(\M{u}^t-\M{u}^\star)\Big\|_{\m{F}}^2.
\end{align}
Hence, combining \eqref{PROOF-IEQ3-NEW2} and \eqref{PROOF-IEQ4}, it holds that
\begin{align}\label{PROOF-IEQ5}
&\mathbb{E} \left[ \left\|\hat{\M{z}}^{t+1}-\M{z}^{\star} \right\|_{\m{F}}^2~|~\theta_t \right]\nonumber\\
&=\Big\|\M{v}^t-\M{v}^\star\Big\|_{\m{F}}^2-\alpha^2\Big\|\sqrt{\nicefrac{1}{2}(\M{I}-\bar{\M{M}})}(\M{u}^t-\M{u}^\star)\Big\|_{\m{F}}^2\nonumber\\
&\quad+\frac{\alpha^2}{p} \Big\|\sqrt{\nicefrac{1}{2}(\M{I}-\bar{\M{M}})}(\hat{\M{u}}^{t+1}-\M{u}^t) \Big\|_{\m{F}}^2\nonumber\\
&\quad-2\alpha \Big\langle\sqrt{\nicefrac{1}{2}(\M{I}-\bar{\M{M}})}(\hat{\M{u}}^{t+1}-\M{u}^{\star}),\M{z}^t-\M{z}^{\star}\Big\rangle.
\end{align}
From the iteration of $\textsc{MG-Skip}$ \eqref{MG-SKIP-NEW-ITER}, we have
\begin{align*}
\sqrt{\nicefrac{1}{2}(\M{I}-\bar{\M{M}})}\M{u}^{t+1}&=\sqrt{\nicefrac{1}{2}(\M{I}-\bar{\M{M}})}\M{u}^{t}+\nicefrac{p}{\alpha} \bar{\M{z}}^{t+1},
\end{align*}
which implies that
$$
\bar{\M{z}}^{t+1}=\frac{\alpha}{p}\sqrt{\nicefrac{1}{2}(\M{I}-\bar{\M{M}})}\left(\M{u}^{t+1}-\M{u}^{t}\right).
$$
Thus, considering that $\M{x}^{t+1}=\mathrm{prox}_{\alpha R}\left(\M{z}^t-\bar{\M{z}}^{t+1}\right)$, we have
$$
\M{x}^{t+1}=\mathrm{prox}_{\alpha R}\left(\M{z}^t-\nicefrac{\alpha}{p}\sqrt{\nicefrac{1}{2}(\M{I}-\bar{\M{M}})}\left(\M{u}^{t+1}-\M{u}^{t}\right)\right).
$$
Since
\begin{align*}
\hat{\M{u}}^{t+1}&=\M{u}^t+\frac{p}{\alpha} \sqrt{\nicefrac{1}{2}(\M{I}-\bar{\M{M}})}{\M{z}}^{t},\\
\M{u}^{t+1}&=\M{u}^t+\frac{p\theta_t}{\alpha} \sqrt{\nicefrac{1}{2}(\M{I}-\bar{\M{M}})}{\M{z}}^{t},
\end{align*}
it gives that
$$
\M{u}^{t+1}-\M{u}^t=\theta_t(\hat{\M{u}}^{t+1}-\M{u}^t).
$$
Considering that we define
$$
\hat{\M{z}}^{t+1}={\M{z}}^t-\frac{\theta_t\alpha}{p}\sqrt{\nicefrac{1}{2}(\M{I}-\bar{\M{M}})}\left(\hat{\M{u}}^{t+1}-\M{u}^t\right),
$$ it holds that
\begin{align*}
\M{x}^{t+1}&=\mathrm{prox}_{\alpha R}\left(\M{z}^t-\frac{\alpha}{p}\sqrt{\nicefrac{1}{2}(\M{I}-\bar{\M{M}})}\left(\M{u}^{t+1}-\M{u}^{t}\right)\right)\\
&=\mathrm{prox}_{\alpha R}\left(\M{z}^t-\frac{\alpha\theta_t}{p}\sqrt{\nicefrac{1}{2}(\M{I}-\bar{\M{M}})}\left(\hat{\M{u}}^{t+1}-\M{u}^{t}\right)\right)\\
&=\m{prox}_{\alpha R}(\hat{\M{z}}^{t+1}).
\end{align*}
Since $\M{x}^{\star}=\m{prox}_{\alpha R}({\M{z}}^{\star})$, by the nonexpansivity of proximal mapping, we have
$$
\mathbb{E}\left[\left\|\M{x}^{t+1}-\M{x}^{\star}\right\|_{\m{F}}^2~|~\theta_t\right]\leq\mathbb{E}\left[\left\|\hat{\M{z}}^{t+1}-\M{z}^{\star}\right\|_{\m{F}}^2~|~\theta_t\right].
$$
Then, from \eqref{PROOF-IEQ5}, it gives that
\begin{align}\label{PROOF-IEQ6}
&\mathbb{E}\left[\left\|\M{x}^{t+1}-\M{x}^{\star}\right\|_{\m{F}}^2\right]\nonumber\\
&\leq\Big\|\M{v}^t-\M{v}^\star\Big\|_{\m{F}}^2
-\alpha^2\Big\|\sqrt{\nicefrac{1}{2}(\M{I}-\bar{\M{M}})}(\M{u}^t-\M{u}^\star)\Big\|_{\m{F}}^2\nonumber\\
&\quad+\frac{\alpha^2}{p} \Big\|\sqrt{\nicefrac{1}{2}(\M{I}-\bar{\M{M}})}(\hat{\M{u}}^{t+1}-\M{u}^t) \Big\|_{\m{F}}^2\nonumber\\
&\quad-2\alpha \Big\langle\sqrt{\nicefrac{1}{2}(\M{I}-\bar{\M{M}})}(\hat{\M{u}}^{t+1}-\M{u}^{\star}),\M{z}^t-\M{z}^{\star}\Big\rangle.
\end{align}
On the other hand, we consider the dual variable. Taking the expectation over $\theta_t$, we have
\begin{align*}
&\mathbb{E} \left[\left\|\M{u}^{t+1}-\M{u}^\star\right\|_{\m{F}}^2 ~|~\theta_t \right]\nonumber\\
&=p\Big\|\hat{\M{u}}^{t+1}-\M{u}^{\star}\Big\|_{\m{F}}^2+(1-p)\Big\|\M{u}^{t}-\M{u}^{\star}\Big\|_{\m{F}}^2.
\end{align*}
By rewriting $\Big\|\hat{\M{u}}^{t+1}-\M{u}^{\star}\Big\|_{\m{F}}^2$, we have
\begin{align*}
&\Big\|\hat{\M{u}}^{t+1}-\M{u}^{\star}\Big\|_{\m{F}}^2\\
&=\Big\|(\M{u}^t-\M{u}^{\star})+(\hat{\M{u}}^{t+1}-\M{u}^{t})\Big\|_{\m{F}}^2\\
&=\Big\|\M{u}^t-\M{u}^{\star}\Big\|_{\m{F}}^2+\Big\|\hat{\M{u}}^{t+1}-\M{u}^{t}\Big\|^2_{\m{F}}+2\Big\langle\M{u}^t-\M{u}^{\star},\hat{\M{u}}^{t+1}-\M{u}^{t}\Big\rangle\\
&=\Big\|\M{u}^t-\M{u}^{\star}\Big\|_{\m{F}}^2+\Big\|\hat{\M{u}}^{t+1}-\M{u}^{t}\Big\|^2_{\m{F}}\nonumber\\
&\quad+2\Big\langle\M{u}^t-\hat{\M{u}}^{t+1},\hat{\M{u}}^{t+1}-\M{u}^{t}\Big\rangle
+2\Big\langle\hat{\M{u}}^{t+1}-\M{u}^{\star},\hat{\M{u}}^{t+1}-\M{u}^{t}\Big\rangle\\
&=\Big\|\M{u}^t-\M{u}^{\star}\Big\|_{\m{F}}^2-\Big\|\hat{\M{u}}^{t+1}-\M{u}^{t}\Big\|_{\m{F}}^2\\
&\quad+2\Big\langle\hat{\M{u}}^{t+1}-\M{u}^{\star},\hat{\M{u}}^{t+1}-\M{u}^{t}\Big\rangle.
\end{align*}
Noting that
$$
\hat{\M{u}}^{t+1}-\M{u}^{t}=\frac{p}{\alpha}\sqrt{\nicefrac{1}{2}(\M{I}-\bar{\M{M}})}\left(\M{z}^t-\M{z}^\star\right),
$$
it gives that
\begin{align*}
&\Big\|\hat{\M{u}}^{t+1}-\M{u}^{\star}\Big\|_{\m{F}}^2\\
&=\left\|\M{u}^t-\M{u}^{\star}\right\|_{\m{F}}^2-\left\|\hat{\M{u}}^{t+1}-\M{u}^{t}\right\|_{\m{F}}^2\\
&\quad+\frac{2p}{\alpha}\left\langle\hat{\M{u}}^{t+1}-\M{u}^{\star}, \sqrt{\nicefrac{1}{2}(\M{I}-\bar{\M{M}})}({\M{z}}^t-\M{z}^{\star})\right\rangle.
\end{align*}
Therefore, from
$
\mathbb{E} [\|\M{u}^{t+1}-\M{u}^\star\|_{\m{F}}^2 ~|~\theta_t ]
=p\|\hat{\M{u}}^{t+1}-\M{u}^{\star}\|_{\m{F}}^2+(1-p)\|\M{u}^{t}-\M{u}^{\star}\|_{\m{F}}^2
$, we have
\begin{align}\label{PROOF-IEQ7}
&\mathbb{E} \left[\left\|\M{u}^{t+1}-\M{u}^\star\right\|_{\m{F}}^2 \right]= \left\|\M{u}^t-\M{u}^{\star} \right\|_{\m{F}}^2-p \left\|\hat{\M{u}}^{t+1}-\M{u}^{t} \right\|_{\m{F}}^2\nonumber\\
&\quad+\frac{2p^2}{\alpha}  \left\langle\hat{\M{u}}^{t+1}-\M{u}^{\star}, \sqrt{\nicefrac{1}{2}(\M{I}-\bar{\M{M}})}({\M{z}}^t-\M{z}^{\star})\right\rangle.
\end{align}
To eliminate the cross term
$$
-2\alpha \Big\langle\sqrt{\nicefrac{1}{2}(\M{I}-\bar{\M{M}})}(\hat{\M{u}}^{t+1}-\M{u}^{\star}),\M{z}^t-\M{z}^{\star}\Big\rangle
$$
in \eqref{PROOF-IEQ6}, we multiply $\frac{\alpha^2}{p^2}$ on both sides of \eqref{PROOF-IEQ7}, and then adding it to \eqref{PROOF-IEQ6}, it gives that
\begin{align*}
&\mathbb{E}\left[\left\|\M{x}^{t+1}-\M{x}^\star~|~\theta_t \right\|_{\m{F}}^2\right]+\nicefrac{\alpha^2}{ p^2}\mathbb{E}\left[\left\|\M{u}^{t+1}-\M{u}^\star ~|~\theta_t\right\|_{\m{F}}^2\right]\\
&\leq\Big\|\M{v}^t-\M{v}^\star\Big\|_{\m{F}}^2-\alpha^2\Big\|\sqrt{\nicefrac{1}{2}(\M{I}-\bar{\M{M}})}(\M{u}^t-\M{u}^\star)\Big\|_{\m{F}}^2\\
&\quad+\frac{\alpha^2}{p} \Big\|\sqrt{\nicefrac{1}{2}(\M{I}-\bar{\M{M}})}(\hat{\M{u}}^{t+1}-\M{u}^t) \Big\|_{\m{F}}^2\nonumber\\
&\quad-2\alpha \Big\langle\sqrt{\nicefrac{1}{2}(\M{I}-\bar{\M{M}})}(\hat{\M{u}}^{t+1}-\M{u}^{\star}),\M{z}^t-\M{z}^{\star}\Big\rangle\nonumber\\
&\quad+2\alpha \Big\langle\sqrt{\nicefrac{1}{2}(\M{I}-\bar{\M{M}})}(\hat{\M{u}}^{t+1}-\M{u}^{\star}),\M{z}^t-\M{z}^{\star}\Big\rangle\nonumber\\
&\quad+\frac{\alpha^2}{p^2}\Big\|\M{u}^t-\M{u}^{\star} \Big\|_{\m{F}}^2-\frac{\alpha^2}{p} \Big\|\hat{\M{u}}^{t+1}-\M{u}^{t} \Big\|_{\m{F}}^2\nonumber\\
&=\Big\|\M{v}^t-\M{v}^\star\Big\|_{\m{F}}^2+\nicefrac{\alpha^2}{p^2}\Big\|\M{u}^t-\M{u}^{\star} \Big\|_{\m{F}}^2\nonumber\\
&\quad-\alpha^2\Big\|\sqrt{\nicefrac{1}{2}(\M{I}-\bar{\M{M}})}(\M{u}^t-\M{u}^\star)\Big\|_{\m{F}}^2\\
&\quad-\frac{\alpha^2}{p}\left(\M{I}-\left((\M{I}-\bar{\M{M}})/2\right)\right)\Big\|\hat{\M{u}}^{t+1}-\M{u}^{t} \Big\|_{\m{F}}^2.
\end{align*}
Since $-\M{I}\preceq\bar{\M{M}}\preceq\M{I}$, we have $\M{I}-\left((\M{I}-\bar{\M{M}})/2\right)=(\M{I}+\bar{\M{M}})/2\succeq0$.
From \eqref{MG-SKIP-NEW-ITER-dual} and $0\in\mathcal{M}(\M{x}^{\star},\M{z}^{\star},\M{u}^{\star})$, we have, for $\M{u}^0=0$, all the columns of $\M{u}^t$ and $\M{u}^\star$ lie in the column space of $\sqrt{\nicefrac{1}{2}(\M{I}-\bar{\M{M}})}$ for any $t\geq0$. Therefore, it gives that $$
\Big\|\sqrt{\nicefrac{1}{2}(\M{I}-\bar{\M{M}})}(\M{u}^t-\M{u}^{\star})\Big\|_{\m{F}}^2\geq \sigma\Big\|\M{u}^t-\M{u}^{\star}\Big\|_{\m{F}}^2,
$$
where $\sigma=\sigma_m((\M{I}-\bar{\M{M}})/2)\geq\frac{1}{5}$. Then, we can deduce that
\begin{align*}
&\mathbb{E}\left[\left\|\M{x}^{t+1}-\M{x}^\star \right\|_{\m{F}}^2~|~\theta_t\right]+\frac{\alpha^2}{p^2}\mathbb{E}\left[\left\|\M{u}^{t+1}-\M{u}^\star \right\|_{\m{F}}^2~|~\theta_t\right]\nonumber\\
&\leq  \Big\|\M{v}^t-\M{v}^{\star} \Big\|_{\m{F}}^2+\Big(\frac{\alpha^2}{ p^2}-\frac{\alpha^2}{5}\Big) \Big\|\M{u}^t-\M{u}^{\star} \Big\|_{\m{F}}^2,
\end{align*}
i.e., \eqref{PROOF-IEQ8} holds.
\end{IEEEproof}

\section{Proof of Theorem \ref{THM-2}}\label{AP:THM-2}
\begin{IEEEproof}
From the cocoercivity of $(\nabla F-\mu \M{I})$ and \cite[Lemma 3.11]{BOOK2015}, it holds that
\begin{align*}
&\Big\langle\nabla F(\M{x}^t)-\nabla F(\M{x}^\star),\M{x}^t-\M{x}^\star\Big\rangle\\
&\geq\frac{L \mu}{L+\mu}\Big\|\M{x}^t-\M{x}^\star\Big\|_{\m{F}}^2+\frac{1}{L+\mu}\Big\|\nabla F(\M{x}^t)-\nabla F(\M{x}^\star)\Big\|_{\m{F}}^2.
\end{align*}
Thus, it gives that
\begin{align}\label{PROOF-THM2-EQ1}
&\Big\|(\M{x}^t-\alpha\nabla F(\M{x}^t))-(\M{x}^\star-\alpha\nabla F(\M{x}^\star))\Big\|_{\m{F}}^2\nonumber\\
&\leq\left(1-\frac{2\alpha L \mu}{L+\mu}\right)\Big\|\M{x}^t-\M{x}^\star\Big\|_{\m{F}}^2\nonumber\\
&\quad+\left(\alpha^2-\frac{2\alpha}{L+\mu}\right)\Big\|\nabla F(\M{x}^t)-\nabla F(\M{x}^\star)\Big\|_{\m{F}}^2.
\end{align}
On one hand, if $\alpha\leq\frac{2}{L+\mu}$, since
$$
\Big\|\nabla F(\M{x}^t)-\nabla F(\M{x}^\star)\Big\|_{\m{F}}\geq\mu\Big\|\M{x}^t-\M{x}^\star\Big\|_{\m{F}}
$$
and by \eqref{PROOF-THM2-EQ1}, it gives that
\begin{align*}
&\Big\|(\M{x}^t-\alpha\nabla F(\M{x}^t))-(\M{x}^\star-\alpha\nabla F(\M{x}^\star))\Big\|_{\m{F}}^2\\
&\leq\Big(1-\frac{2\alpha L \mu}{L+\mu}+\Big(\alpha^2-\frac{2\alpha}{L+\mu}\Big)\mu^2\Big)\Big\|\M{x}^t-\M{x}^\star\Big\|_{\m{F}}^2\\
&=\Big(1-\alpha\mu\Big)^2\Big\|\M{x}^t-\M{x}^\star\Big\|_{\m{F}}^2.
\end{align*}
On the other hand, if $\alpha\geq\frac{2}{L+\mu}$, since
$$
\Big\|\nabla F(\M{x}^t)-\nabla F(\M{x}^\star)\Big\|_{\m{F}}\leq L\Big\|\M{x}^t-\M{x}^\star\Big\|_{\m{F}}
$$ and by \eqref{PROOF-THM2-EQ1}, it gives that
\begin{align*}
&\Big\|(\M{x}^t-\alpha\nabla f(\M{x}^t))-(\M{x}^\star-\alpha\nabla f(\M{x}^\star))\Big\|_{\m{F}}^2\\
&\leq\Big(1-\frac{2\alpha L \mu}{L+\mu}+\Big(\alpha^2-\frac{2\alpha}{L+\mu}\Big)L^2\Big)\Big\|\M{x}^t-\M{x}^\star\Big\|_{\m{F}}^2\\
&=\Big(1-\alpha L\Big)^2\Big\|\M{x}^t-\M{x}^\star\Big\|_{\m{F}}^2.
\end{align*}
Therefore, we have, when $0<\alpha<\frac{2}{L}$ and $\mu>0$,
\begin{align}\label{EQ:Proof-Theorem1:3}
&\Big\|\M{v}^t-\M{v}^{\star} \Big\|_{\m{F}}^2=\Big\|(\M{x}^t-\alpha\nabla F(\M{x}^t))-(\M{x}^\star-\alpha\nabla F(\M{x}^\star))\Big\|_{\m{F}}^2\nonumber\\
&\leq \max\Big\{(1-\alpha
\mu)^2,(1-\alpha L)^2\Big\}\Big\|\M{x}^t-\M{x}^{\star}\Big\|_{\m{F}}^2.
\end{align}
Combining it with \eqref{PROOF-IEQ8}, we have
\begin{align*}
&\mathbb{E}\left[\left\|\M{x}^{t+1}-\M{x}^\star \right\|_{\m{F}}^2~|~\theta_t\right]+\frac{\alpha^2}{p^2}\mathbb{E}\left[\left\|\M{u}^{t+1}-\M{u}^\star \right\|_{\m{F}}^2~|~\theta_t\right]\nonumber\\
&\leq \zeta \Big(\Big\|\M{x}^t-\M{x}^{\star}\Big\|_{\m{F}}^2+\frac{\alpha^2}{ p^2}\Big\|\M{u}^t-\M{u}^{\star} \Big\|_{\m{F}}^2\Big),
\end{align*}
where
$\zeta:=\max\{(1-\alpha\mu)^2,(1-\alpha L)^2,1-\nicefrac{ p^2}{5}\}<1$.
Taking full expectation and unrolling the recurrence, we have
\begin{align*}
&\mathbb{E}\left[\Psi^{t+1}\right]\leq {\zeta}\mathbb{E}\left[\Psi^{t}\right]\leq{\zeta}^{t+1}\Psi^{0}\\
&\Longrightarrow \mathbb{E}\left[\left\|\M{x}^{t+1}-\M{x}^\star\right\|^2\right]\leq{\zeta}^{t+1}\Psi^{0}.
\end{align*}
Since $\lim_{k\rightarrow\infty}\mathbb{E}\left[\Psi^{t}\right]=0$, using classical results on supermartingale convergence \cite[Proposition A.4.5]{BOOK2015}, we have $\lim_{k\rightarrow\infty}\Psi^{t}=0$ almost surely. Therefore, the proof of Theorem \ref{THM-2} is completed.
\end{IEEEproof}

\end{document}